\documentclass[a4paper,latexcad,twoside,11pt]{article}
\usepackage{amssymb} % \mathfrak, \mathbb
\usepackage{graphicx}
\usepackage{epstopdf}
\def\dfrac{\displaystyle\frac}

\begin{document}
\title{ State generatings for Jones and Kauffman-Jones polynomials }
\author{
Liangxia Wan \thanks{\it E-mail  address: $lxwan@bjtu.edu.cn$. }
 \\
  \small\it Department of Mathematics,
Beijing Jiaotong University, Beijing $100044$, China}

\date{}
\maketitle

\noindent{\small {\bf Abstract} A state generating is introduced to determine the Jones polynomial of a link. Formulae for two infinite families of knots are shown by applying this method, the second family of which  are proved to be non-alternating. Moreover, the method is generalized to compute the Jones-Kauffman polynomial of a virtual link. As examples,  formulae for one infinite family of virtual knots are given.

%\noindent MSC(2000): 05C10, 05C30}
\vskip 2mm
\noindent{\it keywords:} Link; virtual link; state generating; embedding presentation

\vskip 5mm
\noindent {\bf $1.$ Introduction}
\vskip 5mm

\noindent

Given a diagram $L$ of a link in $R^3$ (or $S^3$), denote a crossing by a letter, regard $e=(u^r,v^s)$ as an edge if no any other crossings along the line between $u^r$ and $v^s$, then an embedding presentation $L=(V,E)$ with a rotation ${\cal P}=\sum\limits_{u\in V}\sigma_u$ is obtained \cite{Wan}. Here, $V$ is the set of all crossings and $E$ is the set of edges. $\sigma_u$ is an anticlockwise  rotation of edges incident with $u$. If $e$ is an overcrossing at $u$, then  $r=+$ (omitted for brevity), otherwise $r=-$. Throughout this paper a link $L$ (or a virtual link) is always a corresponding embedding presentation, also a marked diagram (or a marked virtual diagram) unless otherwise specified. The link equivalent class $[L]$ is the corresponding link in $R^3$ (or $S^3$) and the virtual link equivalent class $[L]$ is the corresponding virtual link in $S\times I$.

The Jones polynomial is an invariant of $[L]$ which
brought on major advances in knot theory \cite{Jo85}. %Researchers have done beautiful work on the calculations of Jones polynomials until now.
%Subsequently, Freyd, Yetter, Hoste, Lickorish, Millett and Ocneanu introduced a
The Kauffman bracket polynomial of a link was introduced, which is a simple definition to calculate the corresponding Jones polynomial \cite{Ka87}.  %In \cite{Ka87} Kauffman  introduced the bracket polynomial of a link $L$ which is directly related to the Jones polynomial of $[L]$ .
Based on the Kauffman bracket polynomial, several methods were proposed to compute Jones polynomials of links via
Tutte polynomials \cite{Tu54,Ka89} and Bollob${\rm\acute{a}}$s-
Riordan polynomials \cite{BR01,BR02} for some graphs.
 A spanning tree expansion of Jones polynomial was first introduced by constructing a signed graph in \cite{Th87}. This method was extended in \cite{Ka89}. The Jones polynomial of any link equivalent class can also be calculated from the Bollob${\rm\acute{a}}$s-Riordan polynomial
of the ribbon graph via a certain oriented ribbon graph \cite{DFKLS}. In addition, a matrix for calculating the Jones polynomial of a knot equivalent class was given \cite{Zu95}. %Since the computation of the Jones polynomial is \#P %hard,
However, since
determining the Tutte polynomial of a graph is \#P-hard, it is still tough to calculate the Jones polynomial of a link equivalent class [L] \cite{JVW90}, especially a non-alternating link $L$ with a large crossings.

 A virtual link in $S\times I$ and its Kauffman-Jones polynomial were introduced in \cite{Ka99}, which are the generalizations of a link in $R^3$ (or $S^3$) and its Jones polynomial. Similarly, given a virtual diagram $L$ of a virtual link in $S\times I$, denote a crossing by a letter, regard $e=(u^r,v^s)$ as an edge if no any other crossings along the line between $u^r$ and $v^s$, then an embedding presentation $L=(V,E)$ with a rotation ${\cal P}=\sum\limits_{u\in V}\sigma_u$ is obtained \cite{Wan}. Here, $V$ is the set of all crossings and $E$ is the set of edges. $\sigma_u$ is an anticlockwise  rotation of edges incident with $u$. If $u$ is a classical crossing, then $u^+$ (omitted for brevity) and $u^-$ represent an overcrossing and an undercrossing at $u$ respectively, otherwise $u^+$ and $u^-$ represent two occurrences of $u$.  Throughout this paper a virtual link $L$ is always a corresponding embedding presentation, also a marked virtual diagram unless otherwise specified. The virtual link equivalent class $[L]$ is the corresponding virtual link in $S\times I$.

 Correspondingly, approaches for the Jones polynomial of a link equivalent class [L] were extended to compute the Kauffman-Jones polynomial of a virtual link equivalent class. Firstly, the Kauffman-Jones polynomial of a checkerboard
colorable virtual link $L$ can be calculated via the Bollob${\rm\acute{a}}$s-Riordan poly-
nomial of the corresponding ribbon graph \cite{CP07}. Secondly, a relative variant of the other generalization of the Tutte polynomial can be used to compute the Kauffman-Jones polynomials of some virtual links equivalent classes \cite{DH10}. Thirdly, the Kauffman-Jones polynomial of a virtual link equivalent class was computed via the the signed ribbon graph polynomial of its Seifert ribbon graph \cite{CV08}. In fact, the Jones polynomial of a link equivalent class and the Kauffman-Jones polynomial of a virtual link equivalent class can be computed from the signed ribbon graph polynomial of any of their signed ribbon graphs \cite{Ch09}.

This paper introduces a new method called a state generating to calculate Jones polynomials of links based on their bracket polynomials and generalizes this approach to calculate the Kauffman-Jones polynomial of a virtual link.
\vskip 3mm

\setlength{\unitlength}{0.97mm}
\begin{center}
\begin{picture}(100,30)

\qbezier(0,6)(3,21)(6,24)
\qbezier(0,24)(1,20)(2,16)
\qbezier(3,15)(5,10)(5.5,6)

\qbezier(42,24)(39,15)(42,6)
\qbezier(35,24)(38,15)(35,6)

\put(36.5,15){\circle*{1.5}}
\put(40.5,15){\circle*{1.5}}

\qbezier(70,25)(73,10)(76,25)
\qbezier(70,6)(73,21)(76,6)

\put(73,17.5){\circle*{1.5}}
\put(73,13.5){\circle*{1.5}}

\begin{footnotesize}

\put(4,15){{$u$}}
\put(4.5,19){{$e_1$}}
\put(6.5,24){{$x_1^{r_1}$}}

\put(37.5,14){{$u$}}
\put(72,14.5){{$u^-$}}

\put(-2.5,19){{$e_2$}}
\put(-3.5,24){{$x_2^{r_2}$}}

\put(4.5,11){{$e_4$}}
\put(6.5,4){{$x_4^{r_4}$}}

\put(-2.5,11){{$e_3$}}
\put(-4,4){{$x_3^{r_3}$}}

\put(42,24){{$x_1^{r_1}$}}
\put(31.5,24){{$x_2^{r_2}$}}
\put(31,4){{$x_3^{r_3}$}}
\put(41.5,4){{$x_4^{r_4}$}}

\put(77,24){{$x_1^{r_1}$}}
\put(66.5,24){{$x_2^{r_2}$}}
\put(66,4){{$x_3^{r_3}$}}
\put(76.5,4){{$x_4^{r_4}$}}

\put(30,-1){$A$-splitting }
\put(65,-1){$A^{-1}$-splitting}

\put(26,-7){Fig.$0$: Splitting }

\end{footnotesize}
\end{picture}

\end{center}
\vskip 3mm

 Given a link $L$ with $n$ crossings for $n\ge 2$, let $\sigma_u=(e_1,e_2,e_3,e_4)$ be the rotation at $u\in V(L)$ where $e_1=(u,x_1^{r_1}),e_2=(u^-,x_2^{r_2}),e_3=(u,x_3^{r_3}),e_4=(u^-,x_4^{r_4})$, $r_i\in\{+,-\}$ for $1\le i\le 4$. If one replaces passes $x_1^{r_1}ux_3^{r_3},x_2^{r_2}u^-x_4^{r_4}$ with passes $x_1^{r_1}ux_4^{r_4},x_2^{r_2}ux_3^{r_3}$ respectively, then gets a state $A$ at $u$ denoted by $s_u=A$. Otherwise, if one replaces $x_1^{r_1}ux_3^{r_3},x_2^{r_2}u^-x_4^{r_4}$ with $x_1^{r_1}u^-x_2^{r_2},x_3^{r_3}u^-x_4^{r_4}$ respectively, then gets a state $A^{-1}$ at $u$ denoted by $s_u=A^{-1}$ (See Fig.$0$). 
 By assigning one and only one state of $A$ and $A^{-1}$ to each $u\in V(L)$, one obtains a state $s$ of $L$ and a corresponding graph called the {\it state graph} $G(s)$ of $s$ which consists of loops. Two states $s$ and $s'$ of $L$ are distinct if and only if there exists a crossing $u$ such that $s_u\ne s'_u$. Set $S_L$ to be the set of all of states of $L$. It is obvious that $ S_L$ contains $2^n$ elements. Let $c(s),b(s)$ and $l(s)$ denote the number of crossings, $A^{-1}$ and loops in a state $s$ of $L$ respectively. Then $a(s)=c(s)-b(s)$  is the number of state $A$ in $s$. Let $p_i(L)=\sum\limits_{s\in S_L,l(s)=i}A^{a(s)-b(s)}$. Then the Kauffman bracket polynomial is given below
 $$<L>=\sum\limits_{i\ge 1}p_i(L)(-A^2-A^{-2})^{i-1}.$$
   Thus, the Jones polynomial of $[L]$ is deduced as follow
  $$V_L(t)=(-A)^{-3\omega(L)}<L>
  $$
  where $\omega(L)$ is the writhe of $L$ and $t=A^{-4}$. Let $\rho_h(V_L(t))$ and $\rho_l(V_L(t))$ denote the highest and lowest powers of $t$ occurring in $V_L(t)$ respectively. Then the value
  $br(V_L(t))=\rho_h(V_L(t))$- $\rho_l(V_L(t))$ is called the {\it breath} of $V_L(t)$. Obviously,
  it is enough to calculate $p_i(L)$ and $\omega(L)$ in order to obtain $V_L(t)$.

  Similarly, set $L$ to be a virtual link and set $\Gamma L$ to be a set of its classical crossings with $|\Gamma L|=n$ for $n\ge 1$. Let $\sigma_u=(e_1,e_2,e_3,e_4)$ be the rotation at $u\in V(\Gamma L)$ where $e_1=(u,x_1^{r_1}),e_2=(u^-,x_2^{r_2}),e_3=(u,x_3^{r_3}),e_4=(u^-,x_4^{r_4})$, $r_i\in\{+,-\}$ for $1\le i\le 4$. If one replaces passes $x_1^{r_1}ux_3^{r_3},x_2^{r_2}u^-x_4^{r_4}$ with passes $x_1^{r_1}ux_4^{r_4},x_2^{r_2}ux_3^{r_3}$ respectively, then gets a state $A$ at $u$ denoted by $s_u=A$. Otherwise, if one replaces $x_1^{r_1}ux_3^{r_3},x_2^{r_2}u^-x_4^{r_4}$ with $x_1^{r_1}u^-x_2^{r_2},x_3^{r_3}u^-x_4^{r_4}$ respectively, then gets a state $A^{-1}$ at $u$ denoted by $s_u=A^{-1}$ (See Fig.$0$). By assigning one and only one state of states $A$ and $A^{-1}$ to each $u\in V(\Gamma L)$, one obtains a state $s$ of $L$ and then gets the corresponding {\it state graph} which consists of components. Two states $s$ and $s'$ of $L$ are distinct if and only if there exists a crossing  $u\in V(\Gamma L)$ such that $s_u\ne s'_u$. Set $S_L$ to be the set of all of states of $L$. It is obvious that $ S_L$ contains $2^n$ elements. Let $c(s),b(s)$ and $l(s)$ denote the number of classical crossings, $A^{-1}$ and connected components in a state $s$ of $L$ respectively. Then $a(s)=c(s)-b(s)$ is the number of state $A$ in a state $s$. Let $p_i(L)=\sum\limits_{s\in S_L,l(s)=i}A^{a(s)-b(s)}$. So the Kauffman-Jones polynomial $f_L(A)$ for a virtual link
is given below
 $$f_L(A)=(-A)^{-3\omega(L)}\sum\limits_{i\ge 1}p_i(L)(-A^2-A^{-2})^{i-1}$$
    where $\omega(L)$ is the writhe of $L$.
\vskip 3mm

\setlength{\unitlength}{0.97mm}
\begin{center}
\begin{picture}(100,70)

\put(21.5,55){\circle*{1.5}}
\put(24.5,55){\circle*{1.5}}

\put(43,56.5){\circle*{1.5}}
\put(43,53.5){\circle*{1.5}}

\qbezier(0,46)(3,61)(6,64)
\qbezier(0,64)(1,60)(2,56)
\qbezier(3,55)(5,50)(5.5,46)

\qbezier(0,64)(-6,54)(0,46)
\qbezier(6,64)(12,54)(6,46)

\qbezier(26,64)(23,55)(26,46)
\qbezier(20,64)(23,55)(20,46)

\qbezier(20,64)(14,54)(20,46)
\qbezier(26,64)(32,54)(26,46)

\qbezier(40,64)(43,49)(46,64)
\qbezier(40,46)(43,61)(46,46)

\qbezier(40,64)(34,54)(40,46)
\qbezier(46,64)(52,54)(46,46)

\qbezier(60,58)(63,61)(66,64)
\qbezier(66,58)(64,60)(63,60.5)
\qbezier(62.5,61)(61,62)(60,64)

\qbezier(60,52)(63,55)(66,58)
\qbezier(66,52)(64,54)(63,54.5)
\qbezier(62.5,55)(61,56)(60,58)

\qbezier(60,46)(63,49)(66,52)
\qbezier(66,46)(64,48)(63,48.5)
\qbezier(62.5,49)(61,50)(60,52)

\qbezier(60,64)(54,54)(60,46)
\qbezier(66,64)(72,54)(66,46)

\qbezier(80,58)(83,61)(86,64)
\qbezier(86,58)(84,60)(83,60.5)
\qbezier(82.5,61)(81,62)(80,64)

\qbezier(80,52)(83,55)(86,58)
\qbezier(86,52)(84,54)(83,54.5)
\qbezier(82.5,55)(81,56)(80,58)

\qbezier(80,46)(83,49)(86,52)
\qbezier(86,46)(84,48)(83,48.5)
\qbezier(82.5,49)(81,50)(80,52)

\qbezier(80,64)(74,54)(80,46)
\qbezier(86,64)(92,54)(86,46)

\put(80.5,63){\line(1,0){2}}
\put(80.5,63){\line(0,-1){2}}

\qbezier(-20,34)(-17,31)(-20,28)
\qbezier(-20,28)(-17,25)(-20,22)
\qbezier(-20,22)(-17,19)(-20,16)

\qbezier(-14,34)(-17,31)(-14,28)
\qbezier(-14,28)(-17,25)(-14,22)
\qbezier(-14,22)(-17,19)(-14,16)

\qbezier(-20,34)(-26,24)(-20,16)
\qbezier(-14,34)(-8,24)(-14,16)

\put(-18.5,31){\circle*{1.5}}
\put(-15.5,31){\circle*{1.5}}

\put(-18.5,25){\circle*{1.5}}
\put(-15.5,25){\circle*{1.5}}

\put(-18.5,19){\circle*{1.5}}
\put(-15.5,19){\circle*{1.5}}

\put(1.5,31){\circle*{1.5}}
\put(4.5,31){\circle*{1.5}}

\put(3,26.5){\circle*{1.5}}
\put(3,23.5){\circle*{1.5}}

\put(1.5,19){\circle*{1.5}}
\put(4.5,19){\circle*{1.5}}

\qbezier(0,34)(3,31)(0,28)

\qbezier(0,22)(3,19)(0,16)

\qbezier(6,34)(3,31)(6,28)

\qbezier(6,22)(3,19)(6,16)

\qbezier(0,28)(3,25)(6,28)
\qbezier(6,22)(3,25)(0,22)

\qbezier(0,34)(-6,24)(0,16)
\qbezier(6,34)(12,24)(6,16)

\put(21.5,31){\circle*{1.5}}
\put(24.5,31){\circle*{1.5}}

\qbezier(20,34)(23,31)(20,28)
\qbezier(20,28)(23,25)(20,22)

\qbezier(26,34)(23,31)(26,28)
\qbezier(26,28)(23,25)(26,22)

\qbezier(26,22)(23,19)(20,22)
\qbezier(20,16)(23,19)(26,16)

\qbezier(20,34)(14,24)(20,16)
\qbezier(26,34)(32,24)(26,16)

\put(21.5,25){\circle*{1.5}}
\put(24.5,25){\circle*{1.5}}

\put(23,20.5){\circle*{1.5}}
\put(23,17.5){\circle*{1.5}}

\put(41.5,31){\circle*{1.5}}
\put(44.5,31){\circle*{1.5}}

\put(43,26.5){\circle*{1.5}}
\put(43,23.5){\circle*{1.5}}

\put(43,20.5){\circle*{1.5}}
\put(43,17.5){\circle*{1.5}}

\qbezier(40,34)(43,31)(40,28)

\qbezier(46,34)(43,31)(46,28)

\qbezier(40,28)(43,25)(46,28)
\qbezier(46,22)(43,25)(40,22)

\qbezier(46,22)(43,19)(40,22)
\qbezier(40,16)(43,19)(46,16)

\qbezier(40,34)(34,24)(40,16)
\qbezier(46,34)(52,24)(46,16)

\qbezier(60,28)(63,25)(60,22)
\qbezier(60,22)(63,19)(60,16)

\qbezier(66,28)(63,25)(66,22)
\qbezier(66,22)(63,19)(66,16)

\qbezier(60,34)(63,31)(66,34)
\qbezier(60,28)(63,31)(66,28)

\qbezier(60,34)(54,24)(60,16)
\qbezier(66,34)(72,24)(66,16)

\put(63,32.5){\circle*{1.5}}
\put(63,29.5){\circle*{1.5}}

\put(61.5,25){\circle*{1.5}}
\put(64.5,25){\circle*{1.5}}

\put(61.5,19){\circle*{1.5}}
\put(64.5,19){\circle*{1.5}}

\qbezier(80,22)(83,19)(80,16)

\qbezier(86,22)(83,19)(86,16)

\qbezier(80,34)(83,31)(86,34)
\qbezier(80,28)(83,31)(86,28)

\qbezier(80,28)(83,25)(86,28)
\qbezier(80,22)(83,25)(86,22)

\qbezier(80,34)(74,24)(80,16)
\qbezier(86,34)(92,24)(86,16)

\put(83,32.5){\circle*{1.5}}%
\put(83,29.5){\circle*{1.5}}

\put(83,26.5){\circle*{1.5}}
\put(83,23.5){\circle*{1.5}}

\put(81.5,19){\circle*{1.5}}
\put(84.5,19){\circle*{1.5}}

\qbezier(100,28)(103,25)(100,22)

\qbezier(106,28)(103,25)(106,22)

\qbezier(100,34)(103,31)(106,34)
\qbezier(100,28)(103,31)(106,28)

\qbezier(100,22)(103,19)(106,22)
\qbezier(100,16)(103,19)(106,16)

\qbezier(100,34)(94,24)(100,16)
\qbezier(106,34)(112,24)(106,16)

\put(103,32.5){\circle*{1.5}}
\put(103,29.5){\circle*{1.5}}

\put(101.5,25){\circle*{1.5}}
\put(104.5,25){\circle*{1.5}}

\put(103,20.5){\circle*{1.5}}
\put(103,17.5){\circle*{1.5}}

\qbezier(120,34)(123,31)(126,34)
\qbezier(120,28)(123,31)(126,28)

\qbezier(120,22)(123,19)(126,22)
\qbezier(120,16)(123,19)(126,16)

\qbezier(120,28)(123,25)(126,28)
\qbezier(120,22)(123,25)(126,22)

\qbezier(120,34)(114,24)(120,16)
\qbezier(126,34)(132,24)(126,16)

\put(123,32.5){\circle*{1.5}}
\put(123,29.5){\circle*{1.5}}

\put(123,26.5){\circle*{1.5}}
\put(123,23.5){\circle*{1.5}}

\put(123,20.5){\circle*{1.5}}
\put(123,17.5){\circle*{1.5}}

\begin{footnotesize}

\put(-18.5,33){{$x_1$}}

\put(-18.5,27){{$x_2$}}

\put(-18.5,16){{$x_3$}}

\put(1.5,33){{$x_1$}}

\put(2,24.2){{$x_2^-$}}

\put(1.5,16){{$x_3$}}

\put(21.5,33){{$x_1$}}

\put(21.5,27){{$x_2$}}

\put(22,18.2){{$x_3^-$}}

\put(41.5,33){{$x_1$}}
\put(42,24.2){{$x_2^-$}}
\put(42,18.2){{$x_3^-$}}

\put(62,30.2){{$x_1^-$}}

\put(61.5,23){{$x_2$}}

\put(61.5,16){{$x_3$}}

\put(82,30.2){{$x_1^-$}}

\put(82,24.2){{$x_2^-$}}

\put(81.5,16){{$x_3$}}

\put(102,30.2){{$x_1^-$}}

\put(101.5,23){{$x_2$}}

\put(102,18.2){{$x_3^-$}}

\put(122,30.2){{$x_1^-$}}

\put(122,24.2){{$x_2^-$}}

\put(122,18.2){{$x_3^-$}}

\put(3,55){{$x_1$}}

\put(21.5,51){{$x_1$}}

\put(42.5,54.3){{$x_1^-$}}

\put(63.8,60.5){{$x_1$}}
\put(63.8,54.5){{$x_2$}}
\put(63.8,48.5){{$x_3$}}

\put(83.8,60.5){{$x_1$}}
\put(83.8,54.5){{$x_2$}}
\put(83.8,48.5){{$x_3$}}
\put(0.5,41){{$O$}}
\put(20.5,41){{$s_1$}}
\put(40.5,41){{$s_2$}}
\put(60.5,41){{$RT_0$}}
\put(73.5,41){{oriented $RT_0$}}

\put(-21,11){{$s_1(0)$}}
\put(-1,11){{$s_2(0)$}}
\put(19,11){{$s_3(0)$}}
\put(39,11){{$s_4(0)$}}
\put(59,11){{$s_5(0)$}}
\put(79,11){{$s_6(0)$}}
\put(99,11){{$s_7(0)$}}
\put(119,11){{$s_8(0)$}}
\put(10,1){{Fig.1: Jones polynomial of the right handed trefoil $RT_0$}}
\end{footnotesize}
\end{picture}
\end{center}

 Now we introduce a state generating  to calculate the Jones polynomial of a link and the Kauffman-Jones polynomial for a virtual link. In order to calculate the Jones polynomial of a link $L$ (or a Kauffman-Jones polynomial of a virtual link $L$), choose a link $L_1$ (or a virtual link $L_1$ ) with $|V(L_1)|<|V(L)|$ (or $|V(\Gamma L_1)|<|V(\Gamma L)|$) such that each state of $L$ is generated by some state of $L_1$. This method is called a {\it state generating}. If a state $s_1$ of $L_1$ generates a state of $s$ of $L$, then $s_1$ is called the {\it parent} of $s$ denoted by $par(s)$.

  For example, in order to calculate the Jones polynomial of the right handed trefoil $RT_0$, we choose the unknot $O$ shown in Fig.1. $O$ contains two distinct states $s_j$ whose state graphs are $(x_1)(x_1)$ and $(x_1^-x_1^-)$, respectively, for $1\le j\le 2$. It is clear that $p_1(O)=A^{-1}$ and $p_2(O)=A$. The state $s_1$ generates four distinct states $s_j(0)$ of $RT_0$ for $1\le j\le 4$. The state $s_2$ generates four distinct states $s_j(0)$ of $RT_0$ for $5\le j\le 8$ (See Fig.1). We show their loops of state graphs of $s_j(0)$ with loop number in brackets in sequences below  for $1\le j\le 8$
  \vskip 2mm
 \hskip 5mm $(x_1x_3x_2)(x_1x_2x_3)\{2\}$ \hskip 15mm $(x_1x_2^-x_1x_3x_2^-x_3)\{1\}$ \hskip 15mm
  $(x_1x_2x_3^-x_2x_1x_3^-)\{1\}$

  \hskip 5mm $(x_1x_2^-x_1x_3^-)(x_2^-x_3^-)\{2\}$ \hskip 12mm
  $(x_1^-x_3x_2x_1^-x_2x_3)\{1\}$ \hskip 15mm $(x_1^-x_3x_2^-x_3)(x_1^-x_2^-)\{2\}$

  \hskip 5mm $(x_1^-x_3^-)(x_1^-x_2x_3^-x_2)\{2\}$ \hskip 12mm $(x_1^-x_3^-)(x_1^-x_2^-)(x_2^-x_3^-)\{3\}$

  \vskip 2mm

Obviously,
$$\left\{
\begin{array}{ll}
 p_1(RT_0)=A^{2}p_1(O)+2p_2(O)=A+2A=3A,\\
 p_2(RT_0)=2p_1(O)+(A^2+A^{-2})p_2(O)=3A^{-1}+A^{3},\\
 p_3(RT_0)=A^{-2}p_1(O)=A^{-3}.
  \end{array}
  \right.
 $$
 Then $$<RT_0>=3A+(3A^{-1}+A^{3})(-A^2-A^{-2})+A^{-3}(-A^2-A^{-2})^{2}=A^{-7}-A^{-3}-A^{5}.$$
 Since $\omega(RT_0)=3$,
 $$V_{RT_0}(t)=(-A)^{-9}(A^{-7}-A^{-3}-A^{5})=A^{-4}+A^{-12}-A^{-16}=t+t^{3}-t^{4}.$$
\vskip 3mm
Consider $RT_0$. Add $2n$ crossings $y_i$ on $(x_1^-,x_2)$ in sequence and add $2n$ crossings $z_i$ on $(x_1^-,x_3)$ in sequence for $1\le i\le 2n$, delete edges $(x_1,x_2^-),(x_1^-,x_2),(x_1^-,x_3)$,  and then add edges $(x_2,y_{2n}^-)$, $(x_2^-,z_{2n})$, $(x_3,z_{2n}^-)$, $(y_{2k},y_{2k-1}^-),(y_{2k}^-,y_{2k-1})$,
$(z_{2k},z_{2k-1}^-),(z_{2k}^-,z_{2k-1})$,
$(z_{2k},y_{2k+1}^-),(y_{2k},z_{2k-1}^-)$, $(z_{2k}^-,z_{2k+1})$ and $(y_{2k}^-,y_{2k+1})$ for $1\le k\le n$ where $y_{2n+1}^-=x_2^-$, $y_{2n+1}=x_2$ and $z_{2n+1}=x_3$. A type of knots $RT_n$ are obtained for $n\ge 1$, which belong to the first type of knots called {\it $2$-string alternating} knots. $RT_3$ is shown in Fig.2. The Jones polynomials of $RT_n$ are obtained for $n\ge 1$.

\vskip 3mm

\noindent{\bf Theorem $1.1$} {\it For $n\ge 1$
$$V_{RT_n}(t)=\frac{t^{3n}}{\alpha-\bar\alpha}((t+t^3-t^4)(\alpha^{n+1}-\bar\alpha^{n+1})-
(1+t-t^2)(\alpha^{n}-\bar\alpha^{n}))$$
where
$$
\left\{
   \begin{array}{ll}
     \alpha+\bar\alpha=t^{-2}-t^{-1}+2-t+t^2;\\
     \alpha\cdot\bar\alpha=1.
  \end{array}
\right.
$$}
\vskip 3mm
Given a knot $KV_0$ in Fig.3, delete edges $(x_3^-,x_6^-)$, add $2n$ crossings $y_i$ on $(x_2^-,x_3)$ in sequence, $2n$ crossings $z_i$ on $(x_6,x_4^-)$ in sequence for $1\le i\le 2n$, add edges $(x_2^-,y_{1})$, $(x_6^-,y_1^-)$, $(x_6,z_{1})$, $(y_{2k-1}^-,y_{2k})$, $(y_{2k-1},y_{2k}^-)$, $(y_{2k}^-,y_{2k+1})$, $(y_{2k},z_{2k-1})$, $(z_{2k-1}^-,z_{2k})$, $(z_{2k-1},z_{2k}^-)$, $(z_{2k},z_{2k+1}^-)$, $(z_{2k}^-,y_{2k+1})$ for $1\le k\le n-1$ where $y_{2n+1}=x_3$, $y_{2n+1}^-=x_3^-$ and $z_{2n+1}^-=x_4^-$. Then the second type of knots $KV_n$ are given for $n\ge 1$. $KV_1$ is the knot $10_{152}$ \cite{Ro76}. Each $KV_n$ is non-alternating and its Jones polynomial is shown for $n\ge 1$.

\vskip 3mm
\noindent{\bf Theorem $1.2$}
%{\bf Theorem $3.11$}
{\it $KV_n$ are non-alternating knots for $n\ge 1$.}

\vskip 3mm
\noindent{\bf Theorem $1.3$}
%{\bf Theorem $3.12$}
 {\it For $n\ge 1$,
$$V_{KV_n}(A)=A^{(12n+18)}\sum\limits_{i=1}^3g_i(n)$$
where
\begin{eqnarray*}
g_1(n) & = & (A^4+1+A^{-4})A^{-4n-6}+\sum\limits_{i=0}^{n-1}A^{-4i}((\alpha_1^{n-i}-\bar\alpha_1^{n-i})-(2A^4-A^{-4})
 (\alpha_1^{n-1-i}-\bar\alpha_1^{n-1-i}))\\
 & + & (A^{-2}-2A^{-6}+A^{-10})\sum\limits_{i=0}^{n-1}A^{-4i}(\alpha_1^{n-1-i}-\bar\alpha_1^{n-1-i})\\
 & + & (A^{-6}-A^{-10})\sum\limits_{j=0}^{n-1}A^{-4j}(1+A^{8n-8j-4})+\dfrac{A^2-2A^{-2}
 +A^{-6}}{1-A^8}\sum\limits_{j=0}^nA^{-4j}(1-A^{8n-8j}),
\end{eqnarray*}

%\hskip 10mm
$
g_2(n)  =  \dfrac{A^2-A^{-2}+A^{-6}}{A^{4}+1}(1-A^{8n})+\dfrac{A^6+A^{-6}}{A^{4}+1}(A^{8n+4}-1),
$

%\begin{eqnarray*}
$g_3(n) =  \dfrac{(A^2+A^{-2})(A^4-1+A^{-4})}{\alpha_2-\bar\alpha_2}((1-A^{12}+A^6-A^2)
       (\alpha_2^{n+1}-\bar\alpha_2^{n+1})\\$

\hskip 9mm $+ (A^{12}-A^8+A^{4}-A^2)(\alpha_2^{n}-\bar\alpha_2^{n})),
$%\end{eqnarray*}
\vskip 1mm
$
\left\{
   \begin{array}{ll}
     \alpha_1+\bar\alpha_1=A^{8}+2A^{4}+1-2A^{-4};\\
     \alpha_1\cdot\bar\alpha_1=A^{12}+2A^{8}-2-A^{-4}+A^{-8},
  \end{array}
\right.
$
$
\left\{
   \begin{array}{ll}
     \alpha_2+\bar\alpha_2=A^{8}+A^{4}-1-A^{-4};\\
     \alpha_2\cdot\bar\alpha_2=A^{8}-2A^{4}-2A^{-4}-2A^{-8}.
  \end{array}
\right.
$}
\vskip 3mm
Let $x_3$ be a virtual crossing in  $RT_n$ for $n\ge 0$. Then a type of virtual knots $RT'_n$ are obtained. Their Kauffman-Jones polynomials are as follows for $n\ge 1$.
\vskip 3mm
\noindent{\bf Theorem $1.4$}
%{\bf Theorem $4.2$}
{\it For $n\ge 1$
$$f_{RT'_n}(A)=\frac{A^{-12n}}{\alpha-\bar\alpha}((2A^{-4}-A^{-10})(\alpha^{n+1}-\bar\alpha^{n+1})-
(1-A^{-2}+A^{-6}+A^{-8}-A^{-10})(\alpha^{n}-\bar\alpha^{n}))$$
where
$$
\left\{
   \begin{array}{ll}
     \alpha+\bar\alpha=A^8-A^4+2-A^{-4}+A^{-8};\\
     \alpha\cdot\bar\alpha=1.
  \end{array}
\right.
$$
}

\vskip 3mm

This paper is organized as follows. In Section $2$, we use the state generating method introduced in Section $1$ to study the properties of $RT_n$ and then prove Theorem $1.1$ for $n\ge 1$. In Section $3$, we prove Theorems $1.2$ and $1.3$. In Section $4$ we prove Theorem $1.4$ by applying the state generating method for an infinite family of virtual links $RT_n'$ for $n\ge 1$. Finally some open problems are given in Section $5$.

\vskip 5mm
\noindent{\bf $2.$ Jones polynomials of $RT_n$ }

\vskip 5mm
In this section, we divide the set $S(RT_n)$ of all of states for $RT_n$ into four set $S_j(n)$ for $j\in \{\hbox{\rm I},\hbox{\rm II},\hbox{\rm III},\hbox{\rm IV}\}$ for $n\ge 1$. We study the recursive formulae for $S_j(n)$ and then prove Theorem $1.1$.
\vskip 12mm
%\vskip 2mm

\setlength{\unitlength}{0.97mm}
\begin{center}
\begin{picture}(100,50)

\qbezier(-3.5,17)(2,10)(12,16)
\qbezier(5.5,18)(11,19)(12,16)

\qbezier(4,22)(4,15)(4,14)
\qbezier(3.5,18)(-15,18)(-4,22)
\put(-3.5,19){\line(0,1){5}}

\qbezier(-3.5,22)(23,23)(4.5,26)
\qbezier(4,23)(4,27)(4,28)
\put(3.5,26){\line(-1,0){1}}

\put(0,28.5){\circle*{0.7}}
\put(0,26.5){\circle*{0.7}}
\put(0,24.5){\circle*{0.7}}

\qbezier(4,32)(4,30)(4,29.5)
\qbezier(-1.5,28)(-15,28)(-4,32)
\put(-3.5,29){\line(0,1){7}}

\put(-3.5,26){\line(0,1){2}}

\qbezier(-3.5,32)(23,33)(4.5,36)

\qbezier(4,40)(4,36)(4,33)
\qbezier(3.5,36)(-15,36)(-4,40)
\put(-3.5,37){\line(0,1){7}}

\qbezier(-3.5,40)(23,41)(4.5,44)
\qbezier(4,41)(4,45)(4,46)

\qbezier(3.5,44)(-22,48)(-13,20)
\qbezier(-13,20)(-11,10)(3,12)

\qbezier(-3.5,45)(-2,48)(4,46)

\qbezier(46.5,17)(52,10)(62,16)
\qbezier(55.5,18)(61,19)(62,16)

\qbezier(54,22)(54,15)(54,14)
\qbezier(53.5,18)(35,18)(46,22)
\put(46.5,19){\line(0,1){5}}

\qbezier(46.5,22)(73,23)(54.5,26)
\qbezier(54,23)(54,27)(54,28)
\put(53.5,26){\line(-1,0){1}}

\put(50,28.5){\circle*{0.7}}
\put(50,26.5){\circle*{0.7}}
\put(50,24.5){\circle*{0.7}}

\qbezier(54,32)(54,30)(54,29.5)
\qbezier(48.5,28)(35,28)(46,32)
\put(46.5,29){\line(0,1){7}}

\put(46.5,26){\line(0,1){2}}

\qbezier(46.5,32)(73,33)(54.5,36)

\qbezier(54,40)(54,36)(54,33)
\qbezier(53.5,36)(35,36)(46,40)
\put(46.5,37){\line(0,1){7}}

\qbezier(46.5,40)(73,41)(54.5,44)
\qbezier(54,41)(54,45)(54,46)

\qbezier(53.5,44)(28,48)(37,20)
\qbezier(37,20)(39,10)(53,12)

\qbezier(46.5,45)(48,48)(54,46)

\put(36.6,40){\line(0,1){2}}
\put(36.6,40){\line(1,0){2}}

\qbezier(96.5,17)(102,10)(112,16)
\qbezier(105.5,18)(111,19)(112,16)

\qbezier(104,22)(104,15)(104,14)
\qbezier(103.5,18)(85,18)(96,22)
\put(96.5,19){\line(0,1){6}}

\qbezier(96.5,22)(123,23)(104.5,26)
\qbezier(104,23)(104,27)(104,32)

\qbezier(103.5,26)(85,25)(96,32)

\put(96.5,27){\line(0,1){8}}

\qbezier(96.5,32)(123,33)(104.5,36)

\qbezier(104,40)(104,36)(104,33)
\qbezier(103.5,36)(85,36)(96,40)
\put(96.5,37){\line(0,1){7}}

\qbezier(96.5,40)(123,41)(104.5,44)
\qbezier(104,41)(104,45)(104,46)

\qbezier(103.5,44)(78,48)(87,20)
\qbezier(87,20)(89,10)(103,12)

\qbezier(96.5,45)(98,48)(104,46)

\begin{footnotesize}

\put(-7,45){{$x_{1}$}}
\put(-7,41){{$z_{1}$}}
\put(-7,37){{$z_{2}$}}
\put(-7,32.5){{$z_{3}$}}
\put(-7,26.5){{$z_{4}$}}
\put(4.5,45){{$y_{1}$}}
\put(4.5,41.5){{$y_{2}$}}
\put(4.5,37){{$y_{3}$}}
\put(4.5,30.5){{$y_{4}$}}
\put(4.5,27){{$y_{2n-1}$}}
\put(4.5,23.5){{$y_{2n}$}}
\put(-12,23){{$z_{2n-1}$}}
\put(-8,17){{$z_{2n}$}}
\put(4,11){{$x_{3}$}}
\put(4.5,19){{$x_{2}$}}

\put(43,45){{$x_{1}$}}
\put(43,41){{$z_{1}$}}
\put(43,37){{$z_{2}$}}
\put(43,32.5){{$z_{3}$}}
\put(43,26.5){{$z_{4}$}}
\put(54.5,45){{$y_{1}$}}
\put(54.5,41.5){{$y_{2}$}}
\put(54.5,37){{$y_{3}$}}
\put(54.5,30.5){{$y_{4}$}}
\put(54.5,27){{$y_{2n-1}$}}
\put(54.5,23.5){{$y_{2n}$}}
\put(38,23){{$z_{2n-1}$}}
\put(42,17){{$z_{2n}$}}
\put(54,11){{$x_{3}$}}
\put(54.5,19){{$x_{2}$}}

\put(93,45){{$x_{1}$}}
\put(93,41){{$z_{1}$}}
\put(93,37){{$z_{2}$}}
\put(93,32.5){{$z_{3}$}}
\put(93,27){{$z_{4}$}}
\put(104.5,45){{$y_{1}$}}
\put(104.5,41.5){{$y_{2}$}}
\put(104.5,37){{$y_{3}$}}
\put(104.5,30.5){{$y_{4}$}}
\put(104.5,27){{$y_{5}$}}
\put(104.5,23.5){{$y_{6}$}}
\put(93,23){{$z_{5}$}}
\put(93,16.5){{$z_{6}$}}
\put(104,11){{$x_{3}$}}
\put(104.5,19){{$x_{2}$}}

\put(-5,4){{(a) $RT_{n}$ }}
\put(40,4){{(b)  oriented $RT_{n}$ }}
\put(95,4){{(c) $RT_{3}$ }}
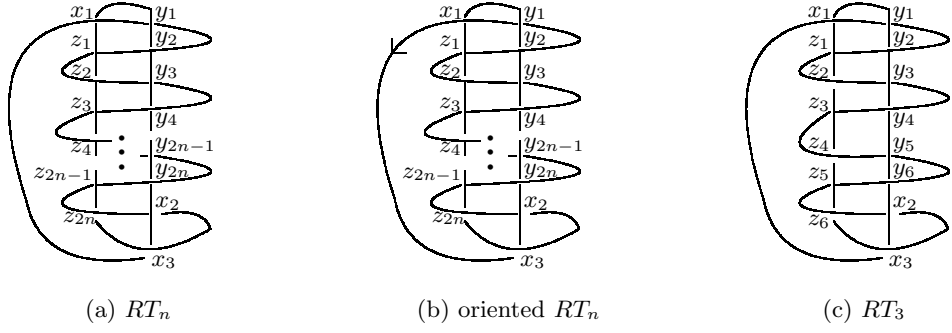
\put(20,-5){{Fig.2: A type of $2$-string alternating knots $RT_{n}$ }}
\end{footnotesize}
\end{picture}
\end{center}
\vskip 3mm

$RT_0$ has eight distinct states $s_j(0)$ shown in Fig.$1$ for $1\le j\le 8$. Each state $s(0)$ of $RT_0$ generates sixteen distinct states of $RT_1$ according to distinct states of $y_i$ and $z_i$ for $1\le i\le 2$. Generally, each state $s(n-1)$ of $RT_{n-1}$ generates sixteen distinct states of $RT_n$ for $n\ge 1$ according to distinct states of $y_i$ and $z_i$ for $2n-1\le i\le 2n$. Let $S({RT}_n)$  denote the set of all of distinct states of $RT_n$ for $n\ge 0$ and set

$$
\begin{array}{ll}
S_{\hbox{\rm\scriptsize I}}(n)=\{s\in S(RT_n),s_{x_2}=s_{x_3}=A|\mbox{ either }\exists 1\le k\le n \mbox{ such that }s_{y_{2k-1}}\cdot s_{y_{2k}}\ne A^2, s_{y_i}=A,\\
 \hskip 18mm s_{z_{2k-1}}=s_{z_{2k}}=s_{z_i}=A^{-}\mbox{ for }2k+1\le i\le 2n\mbox{ or }s_{x_1}=s_{z_i}=A^-, s_{y_i}=A \mbox{ for }1\le i\le 2n \},
\end{array}
$$
$
\begin{array}{ll}
S_{\hbox{\rm\scriptsize II}}(n)=\{s\in S(RT_n),s_{x_2}=s_{x_3}=A|\mbox{ either }\exists 1\le k\le n \mbox{ such that }s_{z_{2k-1}}\cdot s_{z_{2k}}\ne A^{-2},  s_{y_i}=A, \\
\hskip 18mm s_{z_i}=A^{-}
\mbox{ for }2k+1\le i\le 2n\mbox{ or }s_{z_i}=A^-, s_{x_1}=s_{y_i}=A \mbox{ for }1\le i\le 2n \},
\end{array}
$
$$
\begin{array}{ll}
S_{\hbox{\rm\scriptsize III}}(n)=\{s\in S(RT_n),s_{x_2}\cdot s_{x_3}\ne A^2|\mbox{ either }\exists 1\le k\le n \mbox{ such that }s_{y_{2k-1}}\cdot s_{y_{2k}}\ne A^2, s_{y_i}=A,\\
 \hskip 18mm s_{z_{2k-1}}=s_{z_{2k}}=s_{z_i}=A^{-}\mbox{ for }2k+1\le i\le 2n\mbox{ or }s_{x_1}=s_{z_i}=A^-, s_{y_i}=A \mbox{ for }1\le i\le 2n\},
\end{array}
$$
$
\begin{array}{ll}
S_{\hbox{\rm\scriptsize IV}}(n)=\{s\in S(RT_n),s_{x_2}\cdot s_{x_3}\ne A^2|\mbox{ either }\exists 1\le k\le n \mbox{ such that }s_{z_{2k-1}}\cdot s_{z_{2k}}\ne A^{-2},  s_{y_i}=A, \\
\hskip 18mm s_{z_i}=A^{-}
\mbox{ for }2k+1\le i\le 2n\mbox{ or }s_{x_1}=s_{y_i}=A, s_{z_i}=A^-  \mbox{ for }1\le i\le 2n \}.
\end{array}
$

\vskip 2mm

Obviously, there exists one and only one $j\in \{\hbox{\rm I},\hbox{\rm II},\hbox{\rm III},\hbox{\rm IV}\}$  such that $s\in S_j(n)$ for each $s\in S({RT}_n)$. Given $s({n-1})\in S({RT}_{n-1})$, it generates sixteen distinct states $s_j({n})$ of $S({RT}_{n})$ as follows for $1\le j\le 16$:
\vskip 2mm
\hskip 15mm $s_1(n)$ with $s_{y_{2n-1}}=s_{y_{2n}}=s_{z_{2n-1}}=s_{z_{2n}}=A$

\hskip 15mm $s_2(n)$ with $s_{y_{2n-1}}=A^-$ and $s_{y_{2n}}=s_{z_{2n-1}}=s_{z_{2n}}=A$

\hskip 15mm $s_3(n)$ with $s_{y_{2n}}=A^-$ and $s_{y_{2n-1}}=s_{z_{2n-1}}=s_{z_{2n}}=A$

\hskip 15mm $s_4(n)$ with $s_{z_{2n-1}}=A^-$ and $s_{y_{2n-1}}=s_{y_{2n}}=s_{z_{2n}}=A$

\hskip 15mm $s_5(n)$ with $s_{z_{2n}}=A^-$ and $s_{y_{2n-1}}=s_{y_{2n}}=s_{z_{2n-1}}=A$

\hskip 15mm $s_6(n)$ with $s_{y_{2n-1}}=s_{y_{2n}}=A^-$ and $s_{z_{2n-1}}=s_{z_{2n}}=A$

\hskip 15mm $s_7(n)$ with $s_{y_{2n-1}}=s_{z_{2n-1}}=A^-$ and $s_{y_{2n}}=s_{z_{2n}}=A$

\hskip 15mm $s_8(n)$ with $s_{y_{2n-1}}=s_{z_{2n}}=A^-$ and $s_{y_{2n}}=s_{z_{2n-1}}=A$

\hskip 15mm $s_9(n)$ with $s_{y_{2n}}=s_{z_{2n-1}}=A^-$ and $s_{y_{2n-1}}=s_{z_{2n}}=A$

\hskip 15mm $s_{10}(n)$ with $s_{y_{2n}}=s_{z_{2n}}=A^-$ and $s_{y_{2n-1}}=s_{z_{2n-1}}=A$

\hskip 15mm $s_{11}(n)$ with $s_{z_{2n-1}}=s_{z_{2n}}=A^-$ and $s_{y_{2n-1}}=s_{y_{2n}}=A$

\hskip 15mm $s_{12}(n)$ with $s_{y_{2n-1}}=s_{y_{2n}}=s_{z_{2n-1}}=A^-$ and $s_{z_{2n}}=A$

\hskip 15mm $s_{13}(n)$ with $s_{y_{2n-1}}=s_{y_{2n}}=s_{z_{2n}}=A^-$ and $s_{z_{2n-1}}=A$

\hskip 15mm $s_{14}(n)$ with $s_{y_{2n-1}}=s_{z_{2n-1}}=s_{z_{2n}}=A^-$ and $s_{y_{2n}}=A$

\hskip 15mm $s_{15}(n)$ with $s_{y_{2n}}=s_{z_{2n-1}}=s_{z_{2n}}=A^-$ and $s_{y_{2n-1}}=A$

\hskip 15mm $s_{16}(n)$ with $s_{y_{2n-1}}=s_{z_{2n-1}}=s_{z_{2n}}=s_{y_{2n}}=A^-$

\vskip 3mm
\noindent{\bf Lemma $2.1$} {\it Let $s({n-1})\in S_{\hbox{\rm\scriptsize I}}({n-1})$ and let its state graph have $i$ loops for $n,i\ge 1$. Suppose that $s_j(n)$ are states of $RT_n$ and that $par(s_j(n))=s({n-1})$ above for $1\le j\le 16$.

$(1)$ If $j=11$ and $14\le j\le 16$, then $s_j\in S_{\hbox{\rm\scriptsize I}}({n})$. Moreover,
the state graph of $s_{11}(n)$ has $i$ loops, the state graph of $s_{j}(n)$ have $i+1$ loops for $14\le j\le 15$ and the state graph of $s_{16}(n)$ has $i+2$ loops.

$(2)$ Otherwise $s_j(n)\in S_{\hbox{\rm\scriptsize II}}({n})$. Moreover,
the state graph of each $s_{j}(n)$ has $i+1$ loops for $4\le j\le 5$, the state graph of each $s_{j}(n)$ has $i+2$ loops for $j=1$ and $7\le j\le 10$, the state graph of each $s_{16}(n)$ has $i+3$ loops for $2\le j\le 3$ and $12\le j\le 13$ and the state graph of $s_{6}(n)$ has $i+4$ loops.}

\vskip 3mm

{\bf Proof.} Set $s({n-1})\in S_{\hbox{\rm\scriptsize I}}({n-1})$. Without loss of generality, suppose that there exists some $1\le k\le n$ such that $s_{z_{2k-1}}\cdot s_{z_{2k}}=s_{y_{2k-1}}=A^{-}, s_{y_{2k}}= A,\prod\limits_{j=2k+1}^{2n}s_{y_j}=A^{2n-2k},\prod\limits_{j=2k+1}^{2n}s_{z_j}=A^{-(2n-2k)}.$ Other cases are left to readers to verify.

Assume that the state graph of $s(n-1)$ has the loops
$(x_1^{\epsilon}x_3x_2y_{2n-2}Az_{2n-2}^-x_2x_3z_{2n-2}^-B)C$
where $A$ and $B$ are linear sequences, $C$ is the product of $i-1$ loops for $i\ge 1$ and $\epsilon\in\{+,-\}$.

(1) Because $s_{z_{2k-1}}=s_{z_{2k}}=A^-$, $s_{y_{2k-1}}\cdot s_{y_{2k}}\ne A^2$, $s_{z_{i}}=A^-$ and $s_{y_{i}}=A$ for $2k+1\le i\le 2n$ in $s_{11}(n)$,
$$s_{11}(n)\in S_{\hbox{\rm\scriptsize I}}({n}).$$
Because $s_{z_{2n-1}}=s_{z_{2n}}=A^-$ and $s_{y_{2n-1}}\cdot s_{y_{2n}}\ne A^2$ at $s_j(n)$ for $14\le j\le 16$,
$$s_{j}(n)\in S_{\hbox{\rm\scriptsize I}}({n}).$$
Moreover, the state graph of $s_j(n)$ has loops with loop number in brackets in sequence below for $j=11$ and $14\le j\le 16$:
\vskip 2mm
$(x_1^{\epsilon}x_3x_2y_{2n}y_{2n-1}y_{2n-2}Az_{2n-2}^-y_{2n-1}y_{2n}z_{2n-1}^-z_{2n}^-x_2x_3z_{2n}^-z_{2n-1}^-
z_{2n-2}^-B)C\{i\}$%4-11

$(y_{2n-1}^-y_{2n-2}Az_{2n-2}^-)(x_1^{\epsilon}x_3x_2y_{2n}y_{2n-1}^-y_{2n}z_{2n-1}^-z_{2n}^-x_2x_3z_{2n}^-z_{2n-1}^-
z_{2n-2}^-B)C\{i+1\}$%2-7-14

$(y_{2n}^-y_{2n-1}y_{2n-2}Az_{2n-2}^-y_{2n-1})(x_1^{\epsilon}x_3x_2y_{2n}^-z_{2n-1}^-z_{2n}^-x_2x_3z_{2n}^-z_{2n-1}^-
z_{2n-2}^-B)C\{i+1\}$%3-9-15

$(y_{2n}^-y_{2n-1}^-)(y_{2n-2}Az_{2n-2}^-y_{2n-1}^-)(x_1^{\epsilon}x_3x_2y_{2n}^-z_{2n-1}^-z_{2n}^-x_2x_3z_{2n}^-z_{2n-1}^-
z_{2n-2}^-B)C\{i+2\}$%3-9-15-16
\vskip 2mm
Thus the result is clear.

(2) Because $s_{z_{2n-1}}\cdot s_{z_{2n}}\ne A^{-2}$ for $1\le j\le 10$ and $12\le j\le 13$,
$$s_j(n)\in S_{\hbox{\rm\scriptsize II}}({n}).$$
Moreover, the state graph of $s_j(n)$ has loops with loop number in brackets in sequence as follows for $1\le j\le 10$ and $12\le j\le 13$:
\vskip 2mm
$(z_{2n-1}z_{2n})(z_{2n}x_3x_2)(x_1^{\epsilon}x_3x_2y_{2n}y_{2n-1}y_{2n-2}Az_{2n-2}^-y_{2n-1}y_{2n}z_{2n-1}
z_{2n-2}^-B)C\{i+2\}$%1

$(z_{2n-1}z_{2n})(z_{2n}x_3x_2)(y_{2n-1}^-y_{2n-2}Az_{2n-2}^-)(x_1^{\epsilon}x_3x_2y_{2n}y_{2n-1}^-y_{2n}z_{2n-1}
z_{2n-2}^-B)C\{i+3\}$%2

$(z_{2n-1}z_{2n})(z_{2n}x_3x_2)(y_{2n}^-y_{2n-1}y_{2n-2}Az_{2n-2}^-y_{2n-1})(x_1^{\epsilon}x_3x_2y_{2n}^-z_{2n-1}
z_{2n-2}^-B)C\{i+3\}$%3

$(z_{2n}x_3x_2)(x_1^{\epsilon}x_3x_2y_{2n}y_{2n-1}y_{2n-2}Az_{2n-2}^-y_{2n-1}y_{2n}z_{2n-1}^-z_{2n}z_{2n-1}^-
z_{2n-2}^-B)C\{i+1\}$%4

$(z_{2n}^-z_{2n-1}z_{2n}^-x_3x_2)(x_1^{\epsilon}x_3x_2y_{2n}y_{2n-1}y_{2n-2}Az_{2n-2}^-y_{2n-1}y_{2n}z_{2n-1}
z_{2n-2}^-B)C\{i+1\}$%5

$(z_{2n-1}z_{2n})(z_{2n}x_3x_2)(y_{2n-1}^-y_{2n-2}Az_{2n-2}^-)(y_{2n-1}^-y_{2n}^-)(x_1^{\epsilon}x_3x_2y_{2n}^-z_{2n-1}
z_{2n-2}^-B)C\{i+4\}$%2-6

$(z_{2n}x_3x_2)(y_{2n-1}^-y_{2n-2}Az_{2n-2}^-)(x_1^{\epsilon}x_3x_2y_{2n}y_{2n-1}^-y_{2n}z_{2n-1}^-z_{2n}z_{2n-1}^-
z_{2n-2}^-B)C\{i+2\}$%2-7

$(z_{2n}^-z_{2n-1}z_{2n}^-x_3x_2)(y_{2n-1}^-y_{2n-2}Az_{2n-2}^-)(x_1^{\epsilon}x_3x_2y_{2n}y_{2n-1}^-y_{2n}z_{2n-1}
z_{2n-2}^-B)C\{i+2\}$%2-8

$(z_{2n}x_3x_2)(y_{2n}^-y_{2n-1}y_{2n-2}Az_{2n-2}^-y_{2n-1})(x_1^{\epsilon}x_3x_2y_{2n}^-z_{2n-1}^-z_{2n}z_{2n-1}^-
z_{2n-2}^-B)C\{i+2\}$%3-9

$(z_{2n}^-z_{2n-1}z_{2n}^-x_3x_2)(y_{2n}^-y_{2n-1}y_{2n-2}Az_{2n-2}^-y_{2n-1})(x_1^{\epsilon}x_3x_2y_{2n}^-z_{2n-1}
z_{2n-2}^-B)C\{i+2\}$%3-10

$(x_1^{\epsilon}x_3x_2y_{2n}y_{2n-1}y_{2n-2}Az_{2n-2}^-y_{2n-1}y_{2n}z_{2n-1}^-z_{2n}^-x_2x_3z_{2n}^-z_{2n-1}^-
z_{2n-2}^-B)C\{i\}$%4-11

$(z_{2n}x_3x_2)(y_{2n-1}^-y_{2n-2}Az_{2n-2}^-)(y_{2n-1}^-y_{2n}^-)(x_1^{\epsilon}x_3x_2y_{2n}^-z_{2n-1}^-z_{2n}z_{2n-1}^-
z_{2n-2}^-B)C\{i+3\}$%2-6-12

$(z_{2n}^-z_{2n-1}z_{2n}^-x_3x_2)(y_{2n-1}^-y_{2n-2}Az_{2n-2}^-)(y_{2n-1}^-y_{2n}^-)(x_1^{\epsilon}x_3x_2y_{2n}^-z_{2n-1}
z_{2n-2}^-B)C\{i+3\}$%2-6-13
\vskip 2mm
Therefore the consequence holds. $\Box$

\vskip 3mm
\noindent{\bf Lemma $2.2$} {\it Let $s({n-1})\in S_{\hbox{\rm\scriptsize II}}({n-1})$ and its state graph with loops $i$ for $n\ge 1$ and $i\ge 2$. Suppose that $s_j(n)$ are states of $RT_n$ above and that $par(s_j(n))=s({n-1})$ for $1\le j\le 16$.

$(1)$ If $14\le j\le 16$, then $s_j(n)\in S_{\hbox{\rm\scriptsize I}}({n})$. Moreover, the state graph of
$s_{11}(n)$ has $i$ loops, the state graph of each $s_{j}(n)$ have $i-1$ loops for $14\le j\le 15$ and the state graph of $s_{16}(n)$ has $i$ loops.

$(2)$ Otherwise $s_j(n)\in S_{\hbox{\rm\scriptsize II}}({n})$. Moreover,
the state graph of each $s_{j}(n)$ has $i$ loops for $7\le j\le 11$, the state graph of each $s_{j}(n)$ has $i+1$ loops for each $2\le j\le 5$ and $12\le j\le 13$ and the state graph of $s_{j}(n)$ has $i+2$ loops for $j=1,6$.}

\vskip 3mm

{\bf Proof.} Set $s({n-1})\in S_{\hbox{\rm\scriptsize II}}({n-1})$. Without loss of generality, suppose that there exists some $1\le k\le n-1$ such that $s_{z_{2k}}=s_{y_l}=A,s_{z_l}=A^{-}$ for $2k+1\le l\le 2n-2$. Other cases are left to readers to verify.

Assume that the state graph of $s(n-1)$ has the loops
$(x_1^{\epsilon}x_3x_2y_{2n-2}A)(z_{2n-2}^-x_2x_3z_{2n-2}^-B)C$
where $A$ and $B$ are linear sequences, $C$ is the product of $i-2$ loops and $\epsilon\in\{+,-\}$.

(1) Because there exists $n$ such that $s_{z_{2n-1}}=s_{z_{2n}}=A^-$, $s_{y_{2n-1}}\cdot s_{y_{2n}}\ne A^2$ for $14\le j\le 16$,
$$s_{j}(n)\in S_{\hbox{\rm\scriptsize I}}({n}).$$
Moreover, the state graph of $s_j(n)$ has loops with loop number in brackets in sequence below for $14\le j\le 16$:
\vskip 2mm
$
 (z_{2n-2}^-z_{2n-1}^-z_{2n}^-x_3x_2z_{2n}^-z_{2n-1}^-y_{2n}y_{2n-1}^-
 y_{2n}x_2x_3x_1^{\epsilon}A^{re}y_{2n-2}y_{2n-1}^-z_{2n-2}^-B)C\{i-1\}$%2-7-14

 $
(z_{2n-2}^-z_{2n-1}^-z_{2n}^-x_3x_2z_{2n}^-z_{2n-1}^-y_{2n}^-x_2x_3x_1^{\epsilon}
A^{re}y_{2n-2}y_{2n-1}y_{2n}^-y_{2n-1}z_{2n-2}^-B)C\{i-1\}$%3-9-15

$(y_{2n-1}^-y_{2n}^-)
 (z_{2n-2}^-z_{2n-1}^-z_{2n}^-x_3x_2z_{2n}^-z_{2n-1}^-y_{2n}^-x_2x_3x_1^{\epsilon}A^{re}
 y_{2n-2}y_{2n-1}^-z_{2n-2}^-B)C\{i\}$%2-6-12-16

\vskip 2mm
(2) Because $s_{z_{2n-1}}\cdot s_{z_{2n}}\ne A^-$ for $1\le j\le 13$ and $j\ne 11$,
$$s_{j}(n)\in S_{\hbox{\rm\scriptsize II}}({n}).$$
Because there exists $k$ such that $s_{z_{2k}}=s_{y_{l}}=A$ and that $s_{z_{l}}=A^-$ for $2k+1\le l\le 2n$ in $s_{11}(n)$,
$$s_{11}(n)\in S_{\hbox{\rm\scriptsize II}}({n}).$$
Moreover, $s_j(n)$ has loops with loop number in brackets in sequence below for $1\le j\le 13$:
\vskip 2mm
$(z_{2n-1}z_{2n})(z_{2n}x_3x_2)(x_1^{\epsilon}x_3x_2y_{2n}y_{2n-1}y_{2n-2}A)
           (z_{2n-2}^-z_{2n-1}y_{2n}y_{2n-1}z_{2n-2}^-B)C\{i+2\}$%1

$(z_{2n-1}z_{2n})(z_{2n}x_3x_2)
 (z_{2n-2}^-z_{2n-1}y_{2n}y_{2n-1}^-y_{2n}x_2x_3x_1^{\epsilon}A^{re}y_{2n-2}y_{2n-1}^-z_{2n-2}^-B)C\{i+1\}$%2

$(z_{2n-1}z_{2n})(z_{2n}x_3x_2)
(z_{2n-2}^-z_{2n-1}y_{2n}^-x_2x_3x_1^{\epsilon}A^{re}y_{2n-2}y_{2n-1}y_{2n}^-y_{2n-1}z_{2n-2}^-B)C\{i+1\}$%3

$(z_{2n}x_3x_2)(x_1^{\epsilon}x_3x_2y_{2n}y_{2n-1}y_{2n-2}A)
           (z_{2n-2}^-z_{2n-1}^-z_{2n}z_{2n-1}^-y_{2n}y_{2n-1}z_{2n-2}^-B)C\{i+1\}$%4

$(z_{2n}^-z_{2n-1}z_{2n}^-x_3x_2)(x_1^{\epsilon}x_3x_2y_{2n}y_{2n-1}y_{2n-2}A)
           (z_{2n-2}^-z_{2n-1}y_{2n}y_{2n-1}z_{2n-2}^-B)C\{i+1\}$%1-5

$(z_{2n-1}z_{2n})(z_{2n}x_3x_2)(y_{2n-1}^-y_{2n}^-)
 (z_{2n-2}^-z_{2n-1}y_{2n}^-x_2x_3x_1^{\epsilon}A^{re}y_{2n-2}y_{2n-1}^-z_{2n-2}^-B)C\{i+2\}$%2-6

$(z_{2n}x_3x_2)
 (z_{2n-2}^-z_{2n-1}^-z_{2n}z_{2n-1}^-y_{2n}y_{2n-1}^-
 y_{2n}x_2x_3x_1^{\epsilon}A^{re}y_{2n-2}y_{2n-1}^-z_{2n-2}^-B)C\{i\}$%2-7

$(z_{2n}^-z_{2n-1}z_{2n}^-x_3x_2)
 (z_{2n-2}^-z_{2n-1}y_{2n}y_{2n-1}^-y_{2n}x_2x_3x_1^{\epsilon}A^{re}y_{2n-2}y_{2n-1}^-z_{2n-2}^-B)C\{i\}$%2-8

$(z_{2n}x_3x_2)
(z_{2n-2}^-z_{2n-1}^-z_{2n}z_{2n-1}^-y_{2n}^-x_2x_3x_1^{\epsilon}
A^{re}y_{2n-2}y_{2n-1}y_{2n}^-y_{2n-1}z_{2n-2}^-B)C\{i\}$%3-9

$(z_{2n}^-z_{2n-1}z_{2n}^-x_3x_2)
(z_{2n-2}^-z_{2n-1}y_{2n}^-x_2x_3x_1^{\epsilon}A^{re}y_{2n-2}y_{2n-1}y_{2n}^-y_{2n-1}z_{2n-2}^-B)C\{i\}$%3-10

$(x_1^{\epsilon}x_3x_2y_{2n}y_{2n-1}y_{2n-2}A)
           (z_{2n-2}^-z_{2n-1}^-z_{2n}^-x_3x_2z_{2n}^-z_{2n-1}^-y_{2n}y_{2n-1}z_{2n-2}^-B)C\{i\}$%4-11

$(z_{2n}x_3x_2)(y_{2n-1}^-y_{2n}^-)
 (z_{2n-2}^-z_{2n-1}^-z_{2n}z_{2n-1}^-y_{2n}^-x_2x_3x_1^{\epsilon}A^{re}
 y_{2n-2}y_{2n-1}^-z_{2n-2}^-B)C\{i+1\}$%2-6-12

$(z_{2n}^-z_{2n-1}z_{2n}^-x_3x_2)(y_{2n-1}^-y_{2n}^-)
 (z_{2n-2}^-z_{2n-1}y_{2n}^-x_2x_3x_1^{\epsilon}A^{re}y_{2n-2}y_{2n-1}^-z_{2n-2}^-B)C\{i+1\}$%2-6-13
\vskip 2mm
Thus the result is clear. $\Box$

\vskip 3mm
The following results hold by applying a similar way in the argument of the proof of Lemma $2.1$.
\vskip 3mm
\noindent{\bf Lemma $2.3$} {\it Let $s({n-1})\in S_{\hbox{\rm\scriptsize III}}({n-1})$ and let its state graph $i$ loops for $n\ge 1$ and $i\ge 2$. Suppose that $s_j(n)$ are states of $RT_n$ above and that $par(s_j(n))=s({n-1})$ for $1\le j\le 16$.

$(1)$ If $j=11$ and $14\le j\le 16$, then $s_j\in S_{\hbox{\rm\scriptsize III}}({n})$. Moreover,
the state graph of $s_{11}(n)$ has $i$ loops, the state graph of each $s_{j}(n)$ have $i+1$ loops for $14\le j\le 15$ and the state graph of $s_{16}(n)$ has $i+2$ loops.

$(2)$ Otherwise $s_j(n)\in S_{\hbox{\rm\scriptsize IV}}({n})$. Moreover,
the state graph of each $s_{j}(n)$ has $i-1$ loops for $4\le j\le 5$, the state graph of each $s_{j}(n)$ has $i$ loops for $j=1$ and $7\le j\le 10$, the state graph of each $s_{j}(n)$ has $i+1$ loops for $2\le j\le 3$ and $12\le j\le 13$ and the state graph of $s_{6}(n)$ has $i+2$ loops.}

\vskip 3mm
\noindent{\bf Lemma $2.4$} {\it Let $s({n-1})\in S_{\hbox{\rm\scriptsize IV}}({n-1})$ with loops $i$ for $n,i\ge 1$. Suppose that $s_j(n)$ are states of $RT_n$  above and that $par(s_j(n))=s({n-1})$ for $1\le j\le 16$.

$(1)$ If $14\le j\le 16$, then $s_j(n)\in S_{\hbox{\rm\scriptsize III}}({n})$. Moreover,
the state graph of each $s_{j}(n)$ has $i+1$ loops for $14\le j\le 15$ and the state graph of $s_{16}(n)$ has $i+2$ loops.

$(2)$ Otherwise $s_j(n)\in S_{\hbox{\rm\scriptsize IV}}({n})$. Moreover,
the state graph of each $s_{j}(n)$ has $i$ loops for $7\le j\le 11$, the state graph of each $s_{j}(n)$ has $i+1$ loops for $2\le j\le 5$ and $12\le j\le 13$ and the state graph of $s_{j}(n)$ has $i+2$ loops for $j=1,6$.}

\vskip 3mm
%{\bf Proof.} It is shown that By applying the same argument used in Lemma $2.2$. $\Box$
Recursive relations are given below.
\vskip 3mm

\noindent{\bf Lemma $2.5$} {\it Let $p_{1,\hbox{\rm\scriptsize I}}(RT_0)=A$, $p_{2,\hbox{\rm\scriptsize II}}(RT_0)=A^3$, $p_{2,\hbox{\rm\scriptsize III}}(RT_0)=2A^-$, $p_{3,\hbox{\rm\scriptsize III}}(RT_0)=A^{-3}$, $p_{1,\hbox{\rm\scriptsize IV}}(RT_0)=2A$ and $p_{2,\hbox{\rm\scriptsize IV}}(RT_0)=A^-$. Set $p_{i,\hbox{\rm\scriptsize I}}(RT_n)=\sum\limits_{s\in S_{\hbox{\rm\scriptsize I}}(RT_n),l(s)=i}A^{a(s)-b(s)}$,
$p_{i,\hbox{\rm\scriptsize II}}(RT_n)=\sum\limits_{s\in S_{\hbox{\rm\scriptsize II}}(RT_n),l(s)=i}A^{a(s)-b(s)}$,
$p_{i,\hbox{\rm\scriptsize III}}(RT_n)=\sum\limits_{s\in S_{\hbox{\rm\scriptsize III}}(RT_n),l(s)=i}A^{a(s)-b(s)}$,
$p_{i,\hbox{\rm\scriptsize IV}}(RT_n)=\sum\limits_{s\in S_{\hbox{\rm\scriptsize IV}}(RT_n),l(s)=i}A^{a(s)-b(s)}$.
Then for $n\ge 1$
$$
\left\{
\begin{array}{lll}
p_{i,\hbox{\rm\scriptsize I}}(RT_n)=p_{i,\hbox{\rm\scriptsize
      I}}(RT_{n-1})+2A^{-2}p_{i-1,\hbox{\rm\scriptsize I}}(RT_{n-1})+A^{-4}p_{i-2,\hbox{\rm\scriptsize I}}(RT_{n-1})\\
      \hskip 19mm +2A^{-2}p_{i+1,\hbox{\rm\scriptsize II}}(RT_{n-1})+A^{-4}p_{i,\hbox{\rm\scriptsize II}}(RT_{n-1}); \hskip 41mm (1)\\
p_{i,\hbox{\rm\scriptsize II}}(RT_n)=2A^2p_{i-1,\hbox{\rm\scriptsize
      I}}(RT_{n-1})+(A^4+4)p_{i-2,\hbox{\rm\scriptsize I}}(RT_{n-1})\\
      \hskip 19mm +(2A^2+2A^{-2})p_{i-3,\hbox{\rm\scriptsize I}}(RT_{n-1})+p_{i-4,\hbox{\rm\scriptsize I}}(RT_{n-1})+5p_{i,\hbox{\rm\scriptsize II}}(RT_{n-1})\\
      \hskip 19mm +(4A^2+2A^{-2})p_{i-1,\hbox{\rm\scriptsize II}}(RT_{n-1})+(A^4+1)p_{i-2,\hbox{\rm\scriptsize II}}(RT_{n-1}); \hskip 18mm (2)\\
p_{i,\hbox{\rm\scriptsize III}}(RT_n)=p_{i,\hbox{\rm\scriptsize
      III}}(RT_{n-1})+2A^{-2}p_{i-1,\hbox{\rm\scriptsize III}}(RT_{n-1})+A^{-4}p_{i-2,\hbox{\rm\scriptsize III}}(RT_{n-1})\\
      \hskip 19mm +2A^{-2}p_{i-1,\hbox{\rm\scriptsize IV}}(RT_{n-1})+A^{-4}p_{i-2,\hbox{\rm\scriptsize IV}}(RT_{n-1}); \hskip 35.5mm (3)\\
p_{i,\hbox{\rm\scriptsize IV}}(RT_n)=2A^2p_{i+1,\hbox{\rm\scriptsize
      III}}(RT_{n-1})+(A^4+4)p_{i,\hbox{\rm\scriptsize III}}(RT_{n-1})\\
      \hskip 19mm +(2A^2+2A^{-2})p_{i-1,\hbox{\rm\scriptsize III}}(RT_{n-1})+p_{i-2,\hbox{\rm\scriptsize III}}(RT_{n-1})+5p_{i,\hbox{\rm\scriptsize IV}}(RT_{n-1})\\
      \hskip 19mm +(4A^2+2A^{-2})p_{i-1,\hbox{\rm\scriptsize IV}}(RT_{n-1})+(A^4+1)p_{i-2,\hbox{\rm\scriptsize IV}}(RT_{n-1}).\hskip 16.5mm (4)
\end{array}
\right.
$$}
\vskip 3mm

{\bf Proof.} Based on Lemmas
$2.1$ and $2.2$, for $n\geq 1$,
$$
\left\{
\begin{array}{lll}
p_{i,\hbox{\rm\scriptsize I}}(RT_n)=p_{i,\hbox{\rm\scriptsize
      I}}(RT_{n-1})+2A^{-2}p_{i-1,\hbox{\rm\scriptsize I}}(RT_{n-1})+A^{-4}p_{i-2,\hbox{\rm\scriptsize I}}(RT_{n-1})\\
      \hskip 19mm +2A^{-2}p_{i+1,\hbox{\rm\scriptsize II}}(RT_{n-1})+A^{-4}p_{i,\hbox{\rm\scriptsize II}}(RT_{n-1}); \\
p_{i,\hbox{\rm\scriptsize II}}(RT_n)=2A^2p_{i-1,\hbox{\rm\scriptsize
      I}}(RT_{n-1})+(A^4+4)p_{i-2,\hbox{\rm\scriptsize I}}(RT_{n-1})+(2A^2+2A^{-2})p_{i-3,\hbox{\rm\scriptsize I}}(RT_{n-1})\\
      \hskip 19mm +p_{i-4,\hbox{\rm\scriptsize I}}(RT_{n-1})+5p_{i,\hbox{\rm\scriptsize II}}(RT_{n-1})+(4A^2+2A^{-2})p_{i-1,\hbox{\rm\scriptsize II}}(RT_{n-1})\\
      \hskip 19mm +(A^4+1)p_{i-2,\hbox{\rm\scriptsize II}}(RT_{n-1}).
\end{array}
\right.
$$
Similarly, the following result is clear from Lemmas $2.3$ and $2.4$.
$$
\left\{
\begin{array}{lll}
p_{i,\hbox{\rm\scriptsize III}}(RT_n)=p_{i,\hbox{\rm\scriptsize
      III}}(RT_{n-1})+2A^{-2}p_{i-1,\hbox{\rm\scriptsize III}}(RT_{n-1})+A^{-4}p_{i-2,\hbox{\rm\scriptsize III}}(RT_{n-1})\\
      \hskip 19mm +2A^{-2}p_{i-1,\hbox{\rm\scriptsize IV}}(RT_{n-1})+A^{-4}p_{i-2,\hbox{\rm\scriptsize IV}}(RT_{n-1}); \\
p_{i,\hbox{\rm\scriptsize IV}}(RT_n)=2A^2p_{i+1,\hbox{\rm\scriptsize
      III}}(RT_{n-1})+(A^4+4)p_{i,\hbox{\rm\scriptsize III}}(RT_{n-1})+(2A^2+2A^{-2})p_{i-1,\hbox{\rm\scriptsize III}}(RT_{n-1})\\
      \hskip 19mm +p_{i-2,\hbox{\rm\scriptsize III}}(RT_{n-1})+5p_{i,\hbox{\rm\scriptsize IV}}(RT_{n-1})+(4A^2+2A^{-2})p_{i-1,\hbox{\rm\scriptsize IV}}(RT_{n-1})\\
      \hskip 19mm +(A^4+1)p_{i-2,\hbox{\rm\scriptsize IV}}(RT_{n-1}).
\end{array}
\right.
$$
\hskip 150mm $\Box$

\vskip 3mm

{\bf Proof of Theorem $1.1$.} Set $F_1(x,y)=\sum\limits_{i\ge 1,n\ge 0}p_{i,\hbox{\rm\scriptsize I}}(RT_{n})x^{i-1}y^n$, $F_2(x,y)=\sum\limits_{i\ge 1,n\ge 0}p_{i,\hbox{\rm\scriptsize II}}(RT_{n})x^{i-1}y^n$, $F_3(x,y)=\sum\limits_{i\ge 1,n\ge 0}p_{i,\hbox{\rm\scriptsize III}}(RT_{n})x^{i-1}y^n$, $F_4(x,y)=\sum\limits_{i\ge 1,n\ge 0}p_{i,\hbox{\rm\scriptsize IV}}(RT_{n})x^{i-1}y^n$. Let $f_1(x)=\sum\limits_{i\ge 1}p_{i,\hbox{\rm\scriptsize I}}(RT_{n})x^{i-1}$, let $f_2(x)=\sum\limits_{i\ge 1}p_{i,\hbox{\rm\scriptsize II}}(RT_{n})x^{i-1}$, let $f_3(x)=\sum\limits_{i\ge 1}p_{i,\hbox{\rm\scriptsize III}}(RT_{n})x^{i-1}$ and let $f_4(x)=\sum\limits_{i\ge 1}p_{i,\hbox{\rm\scriptsize IV}}(RT_{n})x^{i-1}$. It follows from equations (1-2) that
$$
\left\{
\begin{array}{lll}
\sum\limits_{i\ge 1,n\ge 1}p_{i,\hbox{\rm\scriptsize I}}(RT_n)x^{i}y^n=\sum\limits_{i\ge 1,n\ge 1}p_{i,\hbox{\rm\scriptsize
      I}}(RT_{n-1})x^{i}y^n+2A^{-2}\sum\limits_{i\ge 1,n\ge 1}p_{i-1,\hbox{\rm\scriptsize I}}(RT_{n-1})x^{i}y^n\\
      \hskip 26.5mm +A^{-4}\sum\limits_{i\ge 1,n\ge 1}p_{i-2,\hbox{\rm\scriptsize I}}(RT_{n-1})x^{i}y^n
      +2A^{-2}\sum\limits_{i\ge 1,n\ge 1}p_{i+1,\hbox{\rm\scriptsize II}}(RT_{n-1})x^{i}y^n\\
      \hskip 26.5mm +A^{-4}\sum\limits_{i\ge 1,n\ge 1}p_{i,\hbox{\rm\scriptsize II}}(RT_{n-1})x^{i}y^n; \hskip 73mm (5)\\
\sum\limits_{i\ge 1,n\ge 1}p_{i,\hbox{\rm\scriptsize II}}(RT_n)x^{i-1}y^n=2A^2\sum\limits_{i\ge 1,n\ge 1}p_{i-1,\hbox{\rm\scriptsize
      I}}(RT_{n-1})x^{i-1}y^n
      +(A^4+4)\sum\limits_{i\ge 1,n\ge 1}p_{i-2,\hbox{\rm\scriptsize I}}(RT_{n-1})x^{i-1}y^n \\
      \hskip 27mm +(2A^2+2A^{-2})\sum\limits_{i\ge 1,n\ge 1}p_{i-3,\hbox{\rm\scriptsize I}}(RT_{n-1})x^{i-1}y^n
     +\sum\limits_{i\ge 1,n\ge 1}p_{i-4,\hbox{\rm\scriptsize I}}(RT_{n-1})x^{i-1}y^n\\
     \hskip 27mm +5\sum\limits_{i\ge 1,n\ge 1}p_{i,\hbox{\rm\scriptsize II}}(RT_{n-1})x^{i-1}y^n
     +(4A^2+2A^{-2})\sum\limits_{i\ge 1,n\ge 1}p_{i-1,\hbox{\rm\scriptsize II}}(RT_{n-1})x^{i-1}y^n\\
     \hskip 27mm +(A^4+1)\sum\limits_{i\ge 1,n\ge 1}p_{i-2,\hbox{\rm\scriptsize II}}(RT_{n-1})x^{i-1}y^n. \\
\end{array}
\right.
$$
Since $p_{1,\hbox{\rm\scriptsize I}}(RT_0)=A$ and $p_{2,\hbox{\rm\scriptsize II}}(RT_0)=A^3$,
the set (5) of equations is reduced to the following set of equations
$$
\left\{
\begin{array}{lll}
(xy+2A^{-2}x^2y+A^{-4}x^3y-x)F_1(x,y)+(2A^{-2}y+A^{-4}xy)F_2(x,y)=-Ax;\\
(2A^2xy+(A^4+4)x^2y+(2A^2+2A^{-2})x^3y+x^4y)F_1(x,y)+((5y+4A^2+2A^{-2})xy\\
     \hskip 34.5mm +(A^4+1)x^2y-1)F_2(x,y)=-A^3x.
\end{array}
\right.
$$
Let
$$D=
\left|
\begin{array}{ll}
xy+2A^{-2}x^2y+A^{-4}x^3y-x & 2A^{-2}y+A^{-4}xy\\
2A^2xy+(A^4+4)x^2y+(2A^2+2A^{-2})x^3y+x^4y & (5y+4A^2+2A^{-2})xy+(A^4+1)x^2y-1\\
\end{array}
\right|
$$
\hskip 7.5mm $=-x+(6x-5x^3+x^5)y-xy^2.$

\noindent Then,

$$F_1(x,y)=\frac{1}{D}
%\begin{equation}
\left|
\begin{array}{ll}
-Ax & 2A^{-2}y+A^{-4}xy\\
-A^{3}x & (5y+4A^2+2A^{-2})xy+(A^4+1)x^2y-1\\
\end{array}
\right|
%\end{equation}
$$
$\hskip 42mm =\dfrac{A(1-3y-A^{-2}xy-4A^2xy-(A^4+1)x^2y}{1-(6-5x^2+x^4)y+y^2}$

\noindent and

$$F_2(x,y)=\frac{1}{D}
%\begin{equation}
\left|
\begin{array}{ll}
xy+2A^{-2}x^2y+A^{-4}x^3y-x & -Ax\\
2A^2xy+(A^4+4)x^2y+(2A^2+2A^{-2})x^3y+x^4y & -A^{3}x
\end{array}
\right|
%\end{equation}
$$
$\hskip 42mm =\dfrac{A^2x+A^2xy+(A^4+2)x^2y+(2A^2+A^{-2}x^3)y+x^4y}{1-(6-5x^2+x^4)y+y^2}.$

\noindent Suppose that $1-(6-5x^2+x^4)y+y^2=(1-\alpha y)(1-\bar\alpha y)$ where
$$
\left\{
   \begin{array}{ll}
     \alpha+\bar\alpha=6-5x^2+x^4;\\
     \alpha\cdot\bar\alpha=1.
  \end{array}
\right.
$$
The following equalities can be obtained
\begin{eqnarray*}
F_1(x,y) & = & \dfrac{A(1+(-3-A^{-2}x-4A^2x-(A^4+1)x^2)y}{(1-\alpha y)(1-\bar\alpha y)} \\
      & = & \frac{A}{\alpha-\bar\alpha}((\alpha-3-A^{-2}x-4A^2x-(A^4+1)x^2)\sum\limits_{n\ge 0}\alpha^ny^n\\
      & + & (3+A^{-2}x+4A^2x+(A^4+1)x^2-\bar\alpha)\sum\limits_{n\ge 0}\bar\alpha^ny^n )
\end{eqnarray*}
\noindent and
\begin{eqnarray*}
F_2(x,y) & = & \dfrac{A^3x(1+(1+(A^2+2A^{-2})x+(2+A^{-4})x^2+A^{-2}x^3)y)}{(1-\alpha y)(1-\bar\alpha y)} \\
      & = & \dfrac{A^3x}{\alpha-\bar\alpha}((\alpha+1+(A^2+2A^{-2})x+(2+A^{-4})x^2+A^{-2}x^3)\sum\limits_{n\ge 0}\alpha^ny^n\\
      & + & (-1-(A^2+2A^{-2})x-(2+A^{-4})x^2-A^{-2}x^3-\bar\alpha)\sum\limits_{n\ge 0}\bar\alpha^ny^n ).
\end{eqnarray*}
Thus, for $n\ge 1$
$$f_1(x)=\dfrac{A}{\alpha-\bar\alpha}((\alpha^{n+1}-\bar\alpha^{n+1})+(3+A^{-2}x+4A^2x+(A^4+1)x^2)
(\alpha^{n}-\bar\alpha^{n})) \hskip 15mm (6)$$
and
$$f_2(x)=\dfrac{A^3x}{\alpha-\bar\alpha}((\alpha^{n+1}-\bar\alpha^{n+1})+(1+(A^2+2A^{-2})x+(2+A^{-4})x^2+A^{-2}x^3)
(\alpha^{n}-\bar\alpha^{n})). \hskip 5mm (7)$$
By a similar way, the following equalities can be concluded for $n\ge 1$
$$f_3(x)=\dfrac{2A^-x+A^{-3}x^2}{\alpha-\bar\alpha}((\alpha^{n+1}-\bar\alpha^{n+1})-(3+(4A^2+A^{-2})x+(A^4+1)x^2)
(\alpha^{n}-\bar\alpha^{n})) \hskip 15mm (8)$$
and
$$f_4(x)=\dfrac{2A+A^{-}x}{\alpha-\bar\alpha}((\alpha^{n+1}-\bar\alpha^{n+1})
+(1+(A^2+2A^{-2})x+(2+A^{-4})x^2+A^{-2}x^3)
(\alpha^{n}-\bar\alpha^{n})).\hskip 5mm (9)$$
Since $RT_n$ contains $4n+3$ crossings and $\omega(v)=1$ for each $v\in V(RT_n)$, by setting $x=-A^2-A^{-2}$ and combining with the equalities (6-9), we conclude the following results for $n\ge 1$
\begin{eqnarray*}
V_{RT_n}(t) & = & (-A)^{-(12n+9)}\sum\limits_{j=1}^4f_j(x)\\
          & = & \dfrac{A^{-12n}}{\alpha-\bar\alpha}((A^{-4}+A^{-12}-A^{-16})(\alpha^{n+1}-\bar\alpha^{n+1})
-(1+A^{-4}-A^{-8})
(\alpha^{n}-\bar\alpha^{n}))\\
& = & \dfrac{t^{3n}}{\alpha-\bar\alpha}((t+t^{3}-t^{4})(\alpha^{n+1}-\bar\alpha^{n+1})
-(1+t-t^{2})
(\alpha^{n}-\bar\alpha^{n}))
\end{eqnarray*}
where
$$
\left\{
   \begin{array}{ll}
     \alpha+\bar\alpha=t^{-2}-t^{-1}+2-t+t^2;\\
     \alpha\cdot\bar\alpha=1.
  \end{array}
\right.
$$
\hskip 150mm $\Box$

%For example, the Jones polynomials of $RT_n$ are as follows for $n=1,2,3$:

\vskip 5mm
\noindent{\bf $3.$ Jones polynomials of $KV_n$ }

\vskip 5mm
In this section, for each $KV_n$ with $n\ge 1$, we divide the set $S(KV_n)$ of all of its states into $S^{(j)}(n)$ for $1\le j\le 3$ and obtain some recursive relations. Based on these relations, $KV_n$ is proved to be non-alternating and Theorem 1.2 is concluded.

\vskip 12mm
%\vskip 2mm

\setlength{\unitlength}{0.97mm}
\begin{center}
\begin{picture}(100,35)

\qbezier(-5,33)(-20,33)(-15,19)
\qbezier(-15,25)(-8,30)(7,38)
\qbezier(-17,24)(-28,13)(-8,20)
\qbezier(-8,20)(-1,22)(-3,31.5)
\qbezier(-15,17)(-15.5,14)(-13,12)
\qbezier(2,36)(-23,42)(-26,14)
\qbezier(-26,14)(-9,6)(1,22)
\qbezier(2,22.5)(9,30)(7,38)
\qbezier(-12.5,10.5)(1,0)(2.5,35)

\qbezier(45,33)(26,31)(30,17)
\qbezier(30,16)(29,14)(29.5,12)
\qbezier(30,10)(35,0)(32,20)
\qbezier(32,20)(35,21)(35,18)
\qbezier(35.2,17)(35.5,16)(35.5,15)

\qbezier(30.5,24)(38,30)(57,36)
\qbezier(29,23)(23,15)(32,17)
\qbezier(33,17)(34,17)(36,17.5)

\qbezier(51,35)(27,42)(24,14)

\qbezier(24,14)(23,9)(32,11)
\qbezier(33,11)(34,10.5)(35,11)

\put(41,13){{\circle*{0.5}}}
\put(39.5,12.5){{\circle*{0.5}}}
\put(38,12){{\circle*{0.5}}}

\qbezier(48,32)(48,28)(46,25)
\qbezier(45.5,24)(44,22)(41,19.5)
\qbezier(40.5,19)(39,18)(38,18)

\qbezier(52,34)(55,33)(60,26)
\qbezier(56.5,25)(58,24)(60,26)
\qbezier(55.5,25.5)(53,26)(49,28)
\qbezier(47,28.7)(36,26.7)(58,20.7)
\qbezier(53,18)(60,16)(58,21)
\qbezier(51.5,18.5)(48,19)(44,22)
\qbezier(55,13)(30,22)(43.5,22.5)
\qbezier(48,12)(58,7)(55,13)
\qbezier(47.5,12.5)(47,12.9)(46,13)

\qbezier(58,30)(59,35)(57,36)
\qbezier(57.5,29)(56,25)(54,22)
\qbezier(53.5,21)(52,16)(50,15)
\qbezier(49.5,14.5)(49,13)(47,12)

\qbezier(95,33)(80,33)(85,19)
\qbezier(85,25)(92,30)(107,38)
\qbezier(83,24)(72,13)(90,19.5)
\qbezier(92,20)(99,22)(97,31.5)
\qbezier(85,17)(84.5,14)(87,12)
\qbezier(102,36)(73,42)(74,14)

\qbezier(74,14)(80,6)(90,12.5)
\qbezier(91,13)(98,14)(105,25)
\qbezier(106.5,26.5)(110,33)(107,38)
\qbezier(87.5,10.5)(93,6)(90,22)
\qbezier(90,22)(91,26)(95,22)
\qbezier(95.5,21)(97,20)(99,17.5)
\qbezier(100,17)(111,18)(103,35)

\begin{footnotesize}

\put(-14.5,23){{$x_{1}$}}
\put(-5,33.5){{$x_{2}$}}
\put(-14.5,16){{$x_{3}$}}
\put(-1,37){{$x_{4}$}}
\put(1.5,20.5){{$x_{5}$}}
\put(-16,9.5){{$x_{6}$}}
\put(-10,2){$KV_0$}

\put(52,36){{$x_{4}$}}
\put(45,33.5){{$x_{2}$}}
\put(31,23){{$x_{1}$}}
\put(26.5,15){{$x_{3}$}}
\put(33,14){{$y_{2n}$}}
\put(58.5,29.5){{$x_{5}$}}
\put(49,17){{$z_{2}$}}
\put(48,29.5){{$y_{1}$}}
\put(47,25){{$y_{2}$}}
\put(53,26.5){{$x_{6}$}}
\put(50.5,20.5){{$z_{1}$}}
\put(41,23.5){{$y_{3}$}}
\put(42,19){{$y_{4}$}}
\put(51,14.5){{$z_{3}$}}
\put(46,10){{$z_{4}$}}
\put(33,9){{$z_{2n-1}$}}
\put(26,8){{$z_{2n}$}}
\put(36,16){{$y_{2n-1}$}}
\put(40,2){$KV_n$}

\put(85.5,23){{$x_{1}$}}
\put(95,33.5){{$x_{2}$}}
\put(85.5,16){{$x_{3}$}}
\put(99,37){{$x_{4}$}}
\put(102,25.5){{$x_{5}$}}
\put(98.5,15){{$x_{6}$}}
\put(96.5,21.5){{$y_{1}$}}
\put(91,18){{$y_{2}$}}
\put(90,2){$KV_1$}
\put(25,-5){{Fig.3: The second type of knots $KV_{n}$ }}
\end{footnotesize}
\end{picture}
\end{center}
\vskip 3mm

Let $S(KV_n)$ be the set of all of states of $KV_n$. Denote three sets below
$$
\begin{array}{ll}
S^{(1)}(n)=\{s\in S(KV_n)|s_{x_i}= A^- \mbox{ for }1\le i\le 6\},
\end{array}
$$
$$
\begin{array}{ll}
S^{(2)}(n)=\{s\in S(KV_n)|s_{x_i}= A^- \mbox{ for }1\le i\le 3, \mbox{ and }\prod\limits_{i=4}^{6}s_{x_i}=A^3, \prod\limits_{i=4}^{6}s_{x_i}=A \mbox{ or }\prod\limits_{i=4}^{6}s_{x_i}=A^-\},
\end{array}
$$
$$
\begin{array}{ll}
S^{(3)}(n)=S(KV_n)\setminus \bigcup\limits_{i=1}^2S^{(i)}(n).
\end{array}
$$
Set $p_{i}^{(j)}(n)=\sum\limits_{s\in S^{(j)}(n),l(s)=i}A^{a(s)-b(s)}$ for $1\le j\le 3$ and $i\ge 1$. Obviously,
   $$p_{i}(KV_n)=\sum\limits_{i=1}^3p_{i}^{(j)}(n).$$

\vskip 3mm
{\it Case $1.$} $s\in S^{(1)}(n)$.

Set

$
\begin{array}{ll}
S^{(1)}_{\hbox{\rm\scriptsize I}}(n)=\{s|s_{y_i}=s_{z_i}= A^- \mbox{ for }1\le i\le 2n\},
\end{array}
$

$
\begin{array}{ll}
S^{(1)}_{\hbox{\rm\scriptsize II}}(n)=\{s|\exists 1\le k_0\le k_1\le n \mbox{ such that }s_{z_{2k_0-1}}\cdot s_{z_{2k_0}}\ne  A^{-2}, s_{z_{2k_1-1}}\cdot s_{z_{2k_1}}\ne  A^{-2},\\
\hskip 18mm s_{y_{2k_0-1}}= s_{y_{2k_0}}=s_{y_l}=s_{z_l}= A^{-}, \mbox{ for }1\le l\le 2k_0-2\mbox{ and }2k_1+1\le l\le 2n\},
\end{array}
$

$
\begin{array}{ll}
S^{(1)}_{\hbox{\rm\scriptsize III}}(n)=\{s|\exists 1\le k_0\le k_1\le n \mbox{ such that }s_{y_{2k_0-1}}\cdot s_{y_{2k_0}}\ne  A^{-2}, s_{y_{2k_1-1}}\cdot s_{y_{2k_1}}\ne  A^{-2},\\
\hskip 18mm s_{z_{2k_1-1}}= s_{z_{2k_1}}=s_{y_l}=s_{z_l}= A^{-}, \mbox{ for }1\le l\le 2k_0-2\mbox{ and }2k_1+1\le l\le 2n\},
\end{array}
$

$
\begin{array}{ll}
S^{(1)}_{\hbox{\rm\scriptsize IV}}(n)=\{s|\exists 1\le k_0\le k_1-1\le n \mbox{ such that }s_{z_{2k_0-1}}\cdot s_{z_{2k_0}}\ne  A^{-2}, s_{y_{2k_1-1}}\cdot s_{y_{2k_1}}\ne  A^{-2},\\
\hskip 18mm s_{z_{2k_1-1}}= s_{z_{2k_1}}=s_{y_{2k_0-1}}= s_{y_{2k_0}}=s_{y_l}=s_{z_l}= A^{-}, \mbox{ for }1\le l\le 2k_0-2\mbox{ and }\\
\hskip 18mm 2k_1+1\le l\le 2n\},
\end{array}
$

$
\begin{array}{ll}
S^{(1)}_{\hbox{\rm\scriptsize V}}(n)=\{s|\exists 1\le k_0\le k_1\le n \mbox{ such that }s_{y_{2k_0-1}}\cdot s_{y_{2k_0}}\ne  A^{-2}, s_{z_{2k_1-1}}\cdot s_{z_{2k_1}}\ne  A^{-2},\\
\hskip 18mm s_{y_l}=s_{z_l}= A^{-}, \mbox{ for }1\le l\le 2k_0-2\mbox{ and }2k_1+1\le l\le 2n\}.
\end{array}
$
\vskip 3mm

\vskip 3mm
By a similar way in the argument of the proof  in Lemma 2.5, the following recursive relations are shown.
\vskip 3mm
\noindent{\bf Lemma $3.1$} {\it Let $p_{3,\hbox{\rm\scriptsize I}}^{(1)}(0)=A^{-6}$. Set $p_{i,\hbox{\rm\scriptsize I}}^{(1)}(n)=\sum\limits_{s\in S_{\hbox{\rm\scriptsize I}}^{(1)}(n),l(s)=i}A^{a(s)-b(s)}$,
$p_{i,\hbox{\rm\scriptsize II}}^{(1)}(n)=\sum\limits_{s\in S_{\hbox{\rm\scriptsize II}}^{(1)}(n),l(s)=i}A^{a(s)-b(s)}$,
$p_{i,\hbox{\rm\scriptsize III}}^{(1)}(n)=\sum\limits_{s\in S_{\hbox{\rm\scriptsize III}}^{(1)}(n),l(s)=i}A^{a(s)-b(s)}$,
$p_{i,\hbox{\rm\scriptsize IV}}^{(1)}(n)=\sum\limits_{s\in S_{\hbox{\rm\scriptsize IV}}^{(1)}(n),l(s)=i}A^{a(s)-b(s)}$ and
$p_{i,\hbox{\rm\scriptsize V}}^{(1)}(n)=\sum\limits_{s\in S_{\hbox{\rm\scriptsize V}}^{(1)}(n),l(s)=i}A^{a(s)-b(s)}$. Then
$$p_{i}^{(1)}(n)=p_{i,\hbox{\rm\scriptsize I}}^{(1)}(n)+p_{i,\hbox{\rm\scriptsize II}}^{(1)}(n)
+p_{i,\hbox{\rm\scriptsize III}}^{(1)}(n)+p_{i,\hbox{\rm\scriptsize IV}}^{(1)}(n)+
p_{i,\hbox{\rm\scriptsize V}}^{(1)}(n)$$
where
$$
\left\{
\begin{array}{lll}
p_{i,\hbox{\rm\scriptsize I}}^{(1)}(n)=A^{-4}p_{i,\hbox{\rm\scriptsize
      I}}^{(1)}({n-1}); \\
p_{i,\hbox{\rm\scriptsize II}}^{(1)}(n)=2A^{-2}p_{i+1,\hbox{\rm\scriptsize
      I}}^{(1)}({n-1})+p_{i,\hbox{\rm\scriptsize I}}^{(1)}({n-1})\\
      \hskip 19mm +(A^4+4)p_{i,\hbox{\rm\scriptsize II}}^{(1)}({n-1})+(4A^2+2A^{-2})p_{i-1,\hbox{\rm\scriptsize II}}^{(1)}({n-1})+(A^4+1)p_{i-2,\hbox{\rm\scriptsize II}}^{(1)}({n-1})\\
      \hskip 19mm +2A^{-2}p_{i-1,\hbox{\rm\scriptsize IV}}^{(1)}({n-1})+5p_{i-2,\hbox{\rm\scriptsize IV}}^{(1)}({n-1})\\
      \hskip 19mm +4A^{2}p_{i-3,\hbox{\rm\scriptsize IV}}^{(1)}({n-1})+A^4p_{i-4,\hbox{\rm\scriptsize IV}}^{(1)}({n-1}); \\
p_{i,\hbox{\rm\scriptsize III}}^{(1)}(n)=2A^{-2}p_{i+1,\hbox{\rm\scriptsize
      I}}^{(1)}({n-1})+p_{i,\hbox{\rm\scriptsize I}}^{(1)}({n-1})+A^{-4}p_{i,\hbox{\rm\scriptsize III}}^{(1)}({n-1})\\
      \hskip 19mm +2A^{-2}p_{i-1,\hbox{\rm\scriptsize III}}^{(1)}({n-1})+p_{i-2,\hbox{\rm\scriptsize III}}^{(1)}({n-1})\\
      \hskip 19mm +2A^{-2}p_{i-1,\hbox{\rm\scriptsize
      V}}^{(1)}({n-1})+p_{i-2,\hbox{\rm\scriptsize V}}^{(1)}({n-1}); \\
p_{i,\hbox{\rm\scriptsize IV}}^{(1)}(n)=2A^{-2}p_{i+1,\hbox{\rm\scriptsize
      II}}^{(1)}({n-1})+p_{i,\hbox{\rm\scriptsize II}}^{(1)}({n-1})\\
      \hskip 19mm +A^4p_{i,\hbox{\rm\scriptsize IV}}^{(1)}({n-1})+2A^{-2}p_{i-1,\hbox{\rm\scriptsize
      IV}}^{(1)}({n-1})+p_{i-2,\hbox{\rm\scriptsize IV}}^{(1)}({n-1});\\
p_{i,\hbox{\rm\scriptsize V}}^{(1)}(n)=4p_{i+2,\hbox{\rm\scriptsize
      I}}^{(1)}({n-1})+4A^2p_{i+1,\hbox{\rm\scriptsize I}}^{(1)}({n-1})+A^4p_{i,\hbox{\rm\scriptsize I}}^{(1)}({n-1})+2A^{-2}p_{i+1,\hbox{\rm\scriptsize III}}^{(1)}({n-1})\\
      \hskip 19mm +5p_{i,\hbox{\rm\scriptsize III}}^{(1)}({n-1})+4A^2p_{i-1,\hbox{\rm\scriptsize III}}^{(1)}({n-1})+A^4p_{i-2,\hbox{\rm\scriptsize III}}^{(1)}({n-1})
      +(4+A^{-4})p_{i,\hbox{\rm\scriptsize V}}^{(1)}({n-1})\\
      \hskip 19mm +(4A^2+2A^{-2})p_{i-1,\hbox{\rm\scriptsize V}}^{(1)}({n-1})+(A^4+1)p_{i-2,\hbox{\rm\scriptsize V}}^{(1)}({n-1}).\\
\end{array}
\right.
$$}
\vskip 3mm
{\it Case $2$.} $s\in S^{(2)}(n).$%%%%%%%%%%%%%%%%%%%%%%%%%%%%

Let
\vskip 3mm
 $
\begin{array}{ll}
S^{(2)}_{\hbox{\rm\scriptsize I}}(n)=\{s|\exists 1\le k\le n \mbox{ such that }s_{y_{2k-1}}\cdot s_{y_{2k}}\ne  A^{-2},
 s_{z_{2k-1}}= s_{z_{2k}}=s_{y_l}=s_{z_l}= A^{-},\\
 \hskip 19mm  \mbox{ for }2k+1\le l\le 2n\},
\end{array}
$

$
\begin{array}{ll}
S^{(2)}_{\hbox{\rm\scriptsize II}}(n)=\{s|\mbox{ either }s_{y_i}=s_{z_i}= A^{-}\mbox{ or }\exists 1\le k\le n \mbox{ such that }s_{z_{2k-1}}\cdot s_{z_{2k}}\ne  A^{-2},\\
 \hskip 20mm s_{y_l}=s_{z_l}= A^{-}, \mbox{ for } 1\le i\le 2n, 2k_1+1\le l\le 2n\}.
\end{array}
$

\vskip 3mm
\noindent{\bf Lemma $3.2$} {\it Set $p_{2,\hbox{\rm\scriptsize II}}^{(2)}(0)=3A^{-4}$, $p_{3,\hbox{\rm\scriptsize II}}^{(2)}(0)=3A^{-2}$, $p_{4,\hbox{\rm\scriptsize II}}^{(2)}(0)=1$. Set $p_{i,\hbox{\rm\scriptsize I}}^{(2)}(n)=\sum\limits_{s\in S_{\hbox{\rm\scriptsize I}}^{(2)}(n),l(s)=i}A^{a(s)-b(s)}$,
$p_{i,\hbox{\rm\scriptsize II}}^{(2)}(n)=\sum\limits_{s\in S_{\hbox{\rm\scriptsize II}}^{(2)}(n),l(s)=i}A^{a(s)-b(s)}$. Then
$$p_{i}^{(2)}(n)=p_{i,\hbox{\rm\scriptsize I}}^{(2)}(n)+p_{i,\hbox{\rm\scriptsize II}}^{(2)}(n)$$
where
$$
\left\{
\begin{array}{lll}
p_{i,\hbox{\rm\scriptsize I}}^{(2)}(n)=A^{-4}p_{i,\hbox{\rm\scriptsize I}}^{(2)}({n-1})+2A^{-2}p_{i-1,\hbox{\rm\scriptsize
      I}}^{(2)}({n-1})+p_{i-2,\hbox{\rm\scriptsize I}}^{(2)}({n-1})\\
      \hskip 19mm +2A^{-2}p_{i+1,\hbox{\rm\scriptsize
      II}}^{(2)}({n-1})+p_{i,\hbox{\rm\scriptsize II}}^{(2)}({n-1}); \\
p_{i,\hbox{\rm\scriptsize II}}^{(2)}(n)=2A^{-2}p_{i-1,\hbox{\rm\scriptsize
      I}}^{(2)}({n-1})+5p_{i-2,\hbox{\rm\scriptsize I}}^{(2)}({n-1})+4A^{2}p_{i-3,\hbox{\rm\scriptsize I}}^{(2)}({n-1})\\
      \hskip 19mm +A^{4}p_{i-4,\hbox{\rm\scriptsize I}}^{(2)}({n-1})+(4+A^{-4})p_{i,\hbox{\rm\scriptsize II}}^{(2)}({n-1})\\
      \hskip 19mm +(4A^{2}+2A^{-2})p_{i-1,\hbox{\rm\scriptsize
      II}}^{(2)}({n-1})+(A^4+1)p_{i-2,\hbox{\rm\scriptsize II}}^{(2)}({n-1}).\\
\end{array}
\right.
$$}
\vskip 3mm
{\it Case $3$.} $s\in S^{(3)}(n)$.
\vskip 3mm

Let

 $
\begin{array}{ll}
S^{(3)}_{\hbox{\rm\scriptsize I}}(n)=\{s|\mbox{ either }s_{x_j}=A^-, \prod\limits_{i=1}^3s_{x_i}\ne A^{-3}\mbox{ for }4\le j\le 6\mbox{ or }\exists 1\le k\le n \mbox{ such that }\\
\hskip 20mm s_{y_{2k-1}}\cdot s_{y_{2k}}\ne  A^{-2},
 s_{z_{2k-1}}= s_{z_{2k}}=s_{y_l}=s_{z_l}= A^{-},
 \mbox{ for }2k+1\le l\le 2n\},
\end{array}
$

$
\begin{array}{ll}
S^{(3)}_{\hbox{\rm\scriptsize II}}(n)=\{s|\mbox{ either }s_{y_i}=s_{z_i}= A^{-}, \prod\limits_{i=1}^3s_{x_i}\ne A^{-3}, \prod\limits_{i=4}^6s_{x_i}\ne A^{-3}\mbox{ or }\exists 1\le k\le n \mbox{ such that }\\
 \hskip 20mm s_{z_{2k-1}}\cdot s_{z_{2k}}\ne  A^{-2}, s_{y_l}=s_{z_l}= A^{-}, \mbox{ for } 1\le i\le 2n, 2k+1\le l\le 2n\}.
\end{array}
$
\vskip 3mm
\noindent{\bf Lemma $3.3$} {\it Set $p_{2,\hbox{\rm\scriptsize I}}^{(3)}(0)=3A^{-4}$, $p_{3,\hbox{\rm\scriptsize I}}^{(3)}(0)=3A^{-2}$, $p_{4,\hbox{\rm\scriptsize I}}^{(3)}(0)=1$, $p_{1,\hbox{\rm\scriptsize II}}^{(3)}(0)=9A^{-2}$, $p_{2,\hbox{\rm\scriptsize II}}^{(3)}(0)=18$, $p_{3,\hbox{\rm\scriptsize II}}^{(3)}(0)=15A^2$, $p_{4,\hbox{\rm\scriptsize II}}^{(3)}(0)=6A^4$, $p_{5,\hbox{\rm\scriptsize II}}^{(3)}(0)=A^6$. Set $p_{i,\hbox{\rm\scriptsize I}}^{(3)}(n)=\sum\limits_{s\in S_{\hbox{\rm\scriptsize I}}^{(3)}(n),l(s)=i}A^{a(s)-b(s)}$,
$p_{i,\hbox{\rm\scriptsize II}}^{(3)}(n)=\sum\limits_{s\in S_{\hbox{\rm\scriptsize II}}^{(3)}(n),l(s)=i}A^{a(s)-b(s)}$. Then
$$p_{i}^{(3)}(n)=p_{i,\hbox{\rm\scriptsize I}}^{(3)}(n)+p_{i,\hbox{\rm\scriptsize II}}^{(3)}(n)$$
where
$$
\left\{
\begin{array}{lll}
p_{i,\hbox{\rm\scriptsize I}}^{(3)}(n)=A^{-4}p_{i,\hbox{\rm\scriptsize I}}^{(3)}({n-1})+2A^{-2}p_{i-1,\hbox{\rm\scriptsize
      I}}^{(3)}({n-1})+p_{i-2,\hbox{\rm\scriptsize I}}^{(3)}({n-1})\\
      \hskip 19mm +2A^{-2}p_{i-1,\hbox{\rm\scriptsize
      II}}^{(3)}({n-1})+p_{i-2,\hbox{\rm\scriptsize II}}^{(3)}({n-1}); \\
p_{i,\hbox{\rm\scriptsize II}}^{(3)}(n)=2A^{-2}p_{i+1,\hbox{\rm\scriptsize
      I}}^{(3)}({n-1})+5p_{i,\hbox{\rm\scriptsize I}}^{(3)}({n-1})+4A^{2}p_{i-1,\hbox{\rm\scriptsize I}}^{(3)}({n-1})\\
      \hskip 19mm +A^{4}p_{i-2,\hbox{\rm\scriptsize I}}^{(3)}({n-1})+(4+A^{-4})p_{i,\hbox{\rm\scriptsize II}}^{(2)}({n-1})\\
      \hskip 19mm +(4A^{2}+2A^{-2})p_{i-1,\hbox{\rm\scriptsize
      II}}^{(3)}({n-1})+(A^4+1)p_{i-2,\hbox{\rm\scriptsize II}}^{(3)}({n-1}).\\
\end{array}
\right.
$$}
\hskip 150mm $\Box$

As an example, we calculate the Jones polynomial of $KV_1$ (also $10_{152}$). By Lemma $3.1$, the following results are obtained
$$
\left\{
\begin{array}{lll}
p_{3,\hbox{\rm\scriptsize I}}^{(1)}(1)=A^{-10}, p_{2,\hbox{\rm\scriptsize III}}^{(1)}(1)=p_{2,\hbox{\rm\scriptsize II}}^{(1)}(1)=2A^{-8},
    \\
p_{3,\hbox{\rm\scriptsize III}}^{(1)}(1)=p_{3,\hbox{\rm\scriptsize II}}^{(1)}(1)=A^{-6}, p_{1,\hbox{\rm\scriptsize V}}^{(1)}(1)=4A^{-6}, \hskip 20mm (10)\\
 p_{2,\hbox{\rm\scriptsize V}}^{(1)}(1)=4A^{-4}, p_{3,\hbox{\rm\scriptsize V}}^{(1)}(1)=A^{-2}.\\
\end{array}
\right.
$$
By Lemma $3.2$, the following conclusions are given
$$
\left\{
\begin{array}{lll}
p_{1,\hbox{\rm\scriptsize I}}^{(2)}(1)=6A^{-6}, p_{2,\hbox{\rm\scriptsize I}}^{(2)}(1)=9A^{-4}, p_{3,\hbox{\rm\scriptsize I}}^{(2)}(1)=5A^{-2},
    p_{4,\hbox{\rm\scriptsize I}}^{(2)}(1)=1,
    \\
p_{2,\hbox{\rm\scriptsize II}}^{(2)}(1)=12A^{-4}+3A^{-8}, p_{3,\hbox{\rm\scriptsize II}}^{(2)}(1)=24A^{-2}+9A^{-6}, p_{4,\hbox{\rm\scriptsize II}}^{(2)}(1)=19+10A^{-4}, \hskip 20mm (11)\\
 p_{5,\hbox{\rm\scriptsize II}}^{(2)}(1)=7A^2+5A^{-2}, p_{6,\hbox{\rm\scriptsize II}}^{(2)}(1)=A^{4}+1.\\
\end{array}
\right.
$$
By Lemma $3.3$, we get
$$
\left\{
\begin{array}{lll}
p_{2,\hbox{\rm\scriptsize I}}^{(3)}(1)=18A^{-4}+3A^{-8}, p_{3,\hbox{\rm\scriptsize I}}^{(3)}(1)=9A^{-6}+45A^{-2}, p_{4,\hbox{\rm\scriptsize I}}^{(3)}(1)=48+10A^{-4},
        \\
p_{5,\hbox{\rm\scriptsize I}}^{(3)}(1)=27A^{2}+5A^{-2}, p_{6,\hbox{\rm\scriptsize I}}^{(3)}(1)=8A^{4}+1, p_{7,\hbox{\rm\scriptsize I}}^{(3)}(1)=A^{6}, p_{1,\hbox{\rm\scriptsize II}}^{(3)}(1)=36A^{-2}+15A^{-6},\\ p_{2,\hbox{\rm\scriptsize II}}^{(3)}(1)=108+57A^{-4},
 p_{3,\hbox{\rm\scriptsize II}}^{(3)}(1)=141A^{2}+89A^{-2},p_{4,\hbox{\rm\scriptsize II}}^{(3)}(1)=74+102A^{4}, \hskip 20mm (12)\\
 p_{5,\hbox{\rm\scriptsize II}}^{(3)}(1)=43A^6+35A^{2},p_{6,\hbox{\rm\scriptsize II}}^{(3)}(1)=10A^{8}+9A^{4}, p_{7,\hbox{\rm\scriptsize II}}^{(3)}(1)=A^{10}+A^6.\\
\end{array}
\right.
$$
By combining with the equalities $(10-12)$, we have
$$
\left\{
\begin{array}{lll}
p_{1}(KV_1)=36A^{-2}+25A^{-6},\\
 p_{2}(KV_1)(-A^2-A^{-2})=-108A^2-208A^{-2}-110A^{-6}-10A^{-10}, \\
 p_{3}(KV_1)(-A^2-A^{-2})^2=141A^6+446A^{2}+469A^{-2}+205A^{-6}+22A^{-10}+A^{-14},\\
 p_{4}(KV_1)(-A^2-A^{-2})^3=-102A^{10}-448A^6-752A^2-588A^{-2}-202A^{-6}-20A^{-10},\hskip 7mm (13)\\
 p_{5}(KV_1)(-A^2-A^{-2})^4=43A^{14}+241A^{10}+337A^6+626A^2+379A^{-2}\\
   \hskip 36mm +109A^{-6}+10A^{-10},\\
 p_{6}(KV_1)(-A^2-A^{-2})^5=-10A^{18}-68A^{14}-192A^{10}-290A^6-250A^2\\
   \hskip 36mm -120A^{-2}-28A^{-6}-2A^{-10},\\
 p_{7}(KV_1)(-A^2-A^{-2})^6=A^{22}+8A^{18}+27A^{14}+50A^{10}+55A^6+36A^2+14A^{-2}+2A^{-6}
\end{array}
\right.
$$

By applying the equalities (1) and (13), the Kauffman bracket polynomial of $KV_1$ is
\begin{eqnarray*}
<KV_1> & = & \sum\limits_{i=1}^7p_i(KV_1)(-A^2-A^{-2})^{i-1}\\
     %  & = & \sum\limits_{i=1}^7\sum\limits_{j=1}^3p_i^{(j)}(1)(-A^2-A^{-2})^{i-1}\\
       & = & A^{22}-2A^{18}+2A^{14}-3A^{10}+2A^6-2A^2+A^{-2}+A^{-6}+A^{-14}.
\end{eqnarray*}
Since $\omega(KV_1)=-10$, the Jones polynomial of $RV_1$ is deduced
\begin{eqnarray*}
V_{KV_1}(t) & = & A^{30}(A^{22}-2A^{18}+2A^{14}-3A^{10}+2A^6-2A^2+A^{-2}+A^{-6}+A^{-14})\\
       & = & A^{52}-2A^{48}+2A^{44}-3A^{40}+2A^{36}-2A^{32}+A^{28}+A^{24}+A^{16}\\
       & = & t^{-13}-2t^{-12}+2t^{-11}-3t^{-10}+2t^{-9}-2t^{-8}+t^{-7}+t^{-6}+t^{-4}.
\end{eqnarray*}
\vskip 3mm
It is obvious that $RT_1$ is non-alternating. In order to prove that $RT_n$ are non-alternating for $n\ge 2$. We consider the highest power and the lowest power of $A$ in $f_{i}(n)=p_i(KV_n)(-A^2-A^{-2})^{i-1}$ for $i\ge 1$.
\vskip 3mm
\noindent{\bf Lemma $3.4$} {\it Set $f_{i}(n)=p_i(KV_n)(-A^2-A^{-2})^{i-1}$ for $n,i\ge 1$. Let $\rho_h(f_{n,i})$ and $\rho_l(f_{i}(n))$ denote the highest power and the lowest power of $A$ in $f_{i}(n)$ respectively. Then for $n,i\ge 1$
$$\rho_h(f_{i}(n))\le 8k+14$$
and
$$\rho_l(f_{i}(n))\ge -4k-10.$$}
\vskip 3mm
{\bf Proof.} This conclusion will be verified by induction on $n$. By equalities of $(13)$, the result is obvious for $n=1$.

Assume that the result holds for the integer $k (k\ge 2)$. That is for $i\ge 1$
$$\rho_h(f_{i}(k))\le 8n+14 \mbox{ and }\rho_l(f_{i}(k))\ge -4n-10.$$
This implies
for $i\ge 1$
$$\rho_h(f_{i,j}^{(r)}(k))\le 8n+14 \mbox{ and }\rho_l(f_{i,j}^{(r)}(k))\ge -4n-10 \hskip 10mm (14)$$
where $f_{i,j}^{(r)}(k)=p_i^{(r)}(KV_k)(-A^2-A^{-2})^{i-1}$ for $j\in \{\hbox{\rm I},\hbox{\rm II},\hbox{\rm III}, \hbox{\rm IV}\}$ and $1\le r\le 3$.
By Lemma $3.6$,
$$
\left\{
\begin{array}{lll}
p_{i,\hbox{\rm\scriptsize I}}^{(1)}(k+1)(-A^2-A^{-2})^{i-1}=A^{-4}p_{i,\hbox{\rm\scriptsize
      I}}^{(1)}({k}) (-A^2-A^{-2})^{i-1}; \\
p_{i,\hbox{\rm\scriptsize II}}^{(1)}(k+1)(-A^2-A^{-2})^{i-1}=(2A^{-2}p_{i+1,\hbox{\rm\scriptsize
      I}}^{(1)}({k})+p_{i,\hbox{\rm\scriptsize I}}^{(1)}({k})\\
      \hskip 19mm +(A^4+4)p_{i,\hbox{\rm\scriptsize II}}^{(1)}({k})+(4A^2+2A^{-2})p_{i-1,\hbox{\rm\scriptsize II}}^{(1)}({k})+(A^4+1)p_{i-2,\hbox{\rm\scriptsize II}}^{(1)}({k})\\
      \hskip 19mm +2A^{-2}p_{i-1,\hbox{\rm\scriptsize IV}}^{(1)}({k})+5p_{i-2,\hbox{\rm\scriptsize IV}}^{(1)}({k})\\
      \hskip 19mm +4A^{2}p_{i-3,\hbox{\rm\scriptsize IV}}^{(1)}({k})+A^4p_{i-4,\hbox{\rm\scriptsize IV}}^{(1)}({k}))(-A^2-A^{-2})^{i-1}; \\
p_{i,\hbox{\rm\scriptsize III}}^{(1)}(k+1)(-A^2-A^{-2})^{i-1}=(2A^{-2}p_{i+1,\hbox{\rm\scriptsize
      I}}^{(1)}({k})+p_{i,\hbox{\rm\scriptsize I}}^{(1)}({k})+A^{-4}p_{i,\hbox{\rm\scriptsize III}}^{(1)}({k})\\
      \hskip 19mm +2A^{-2}p_{i-1,\hbox{\rm\scriptsize III}}^{(1)}({k})+p_{i-2,\hbox{\rm\scriptsize III}}^{(1)}({k})\\
      \hskip 19mm +2A^{-2}p_{i-1,\hbox{\rm\scriptsize
      V}}^{(1)}({k})+p_{i-2,\hbox{\rm\scriptsize V}}^{(1)}({k}))(-A^2-A^{-2})^{i-1}; \\
p_{i,\hbox{\rm\scriptsize IV}}^{(1)}(k+1)(-A^2-A^{-2})^{i-1}=(2A^{-2}p_{i+1,\hbox{\rm\scriptsize
      II}}^{(1)}({k})+p_{i,\hbox{\rm\scriptsize II}}^{(1)}({k})\\
      \hskip 19mm +A^4p_{i,\hbox{\rm\scriptsize IV}}^{(1)}({k})+2A^{-2}p_{i-1,\hbox{\rm\scriptsize
      IV}}^{(1)}({k})+p_{i-2,\hbox{\rm\scriptsize IV}}^{(1)}({k}))(-A^2-A^{-2})^{i-1};\\
p_{i,\hbox{\rm\scriptsize V}}^{(1)}(k+1)(-A^2-A^{-2})^{i-1}=(4p_{i+2,\hbox{\rm\scriptsize
      I}}^{(1)}({k})+4A^2p_{i+1,\hbox{\rm\scriptsize I}}^{(1)}({k})+A^4p_{i,\hbox{\rm\scriptsize I}}^{(1)}({k})+2A^{-2}p_{i+1,\hbox{\rm\scriptsize III}}^{(1)}({k})\\
      \hskip 19mm +5p_{i,\hbox{\rm\scriptsize III}}^{(1)}({k})+4A^2p_{i-1,\hbox{\rm\scriptsize III}}^{(1)}({k})+A^4p_{i-2,\hbox{\rm\scriptsize III}}^{(1)}({k})\\
      \hskip 19mm +(4+A^{-4})p_{i,\hbox{\rm\scriptsize V}}^{(1)}({k})\\
      \hskip 19mm +(4A^2+2A^{-2})p_{i-1,\hbox{\rm\scriptsize V}}^{(1)}({k})+(A^4+1)p_{i-2,\hbox{\rm\scriptsize V}}^{(1)}({k}))(-A^2-A^{-2})^{i-1}.\\
\end{array}
\right.
$$
This implies the following equalities
$$
\left\{
\begin{array}{lll}
f_{i,\hbox{\rm\scriptsize I}}^{(1)}(k+1)=A^{-4}f_{i,\hbox{\rm\scriptsize
      I}}^{(1)}({k}); \\
f_{i,\hbox{\rm\scriptsize II}}^{(1)}(k+1)=(2A^{-2}f_{i+1,\hbox{\rm\scriptsize
      I}}^{(1)}({k})(-A^2-A^{-2})^{-1}+f_{i,\hbox{\rm\scriptsize I}}^{(1)}({k})+(A^4+4)f_{i,\hbox{\rm\scriptsize II}}^{(1)}({k})\\
      \hskip 19mm +(4A^2+2A^{-2})f_{i-1,\hbox{\rm\scriptsize II}}^{(1)}({k})(-A^2-A^{-2})+(A^4+1)f_{i-2,\hbox{\rm\scriptsize II}}^{(1)}({k})(-A^2-A^{-2})^{2}\\
      \hskip 19mm +2A^{-2}f_{i-1,\hbox{\rm\scriptsize IV}}^{(1)}({k})(-A^2-A^{-2})+5f_{i-2,\hbox{\rm\scriptsize IV}}^{(1)}({k})(-A^2-A^{-2})^{2}\\
      \hskip 19mm +4A^{2}f_{i-3,\hbox{\rm\scriptsize IV}}^{(1)}({k})(-A^2-A^{-2})^{3}+A^4f_{i-4,\hbox{\rm\scriptsize IV}}^{(1)}({k})(-A^2-A^{-2})^{4}; \\
f_{i,\hbox{\rm\scriptsize III}}^{(1)}(k+1)=2A^{-2}f_{i+1,\hbox{\rm\scriptsize
      I}}^{(1)}({k})(-A^2-A^{-2})^{-1}+f_{i,\hbox{\rm\scriptsize I}}^{(1)}({k})+A^{-4}f_{i,\hbox{\rm\scriptsize III}}^{(1)}({k})\hskip 30mm (15)\\
      \hskip 19mm +2A^{-2}f_{i-1,\hbox{\rm\scriptsize III}}^{(1)}({k})(-A^2-A^{-2})+f_{i-2,\hbox{\rm\scriptsize III}}^{(1)}({k})(-A^2-A^{-2})^{2}\\
      \hskip 19mm +2A^{-2}f_{i-1,\hbox{\rm\scriptsize
      V}}^{(1)}({k})(-A^2-A^{-2})+f_{i-2,\hbox{\rm\scriptsize V}}^{(1)}({k})(-A^2-A^{-2})^{2}; \\
f_{i,\hbox{\rm\scriptsize IV}}^{(1)}(k+1)=2A^{-2}f_{i+1,\hbox{\rm\scriptsize
      II}}^{(1)}({k})(-A^2-A^{-2})^{-1}+f_{i,\hbox{\rm\scriptsize II}}^{(1)}({k})\\
      \hskip 19mm +A^4f_{i,\hbox{\rm\scriptsize IV}}^{(1)}({k})+2A^{-2}f_{i-1,\hbox{\rm\scriptsize
      IV}}^{(1)}({k})(-A^2-A^{-2})+f_{i-2,\hbox{\rm\scriptsize IV}}^{(1)}({k})(-A^2-A^{-2})^{2};\\
f_{i,\hbox{\rm\scriptsize V}}^{(1)}(k+1)=4f_{i+2,\hbox{\rm\scriptsize
      I}}^{(1)}({k})(-A^2-A^{-2})^{-2}+4A^2f_{i+1,\hbox{\rm\scriptsize I}}^{(1)}({k})(-A^2-A^{-2})^{-1}+A^4f_{i,\hbox{\rm\scriptsize I}}^{(1)}({k})\\
      \hskip 19mm +2A^{-2}f_{i+1,\hbox{\rm\scriptsize III}}^{(1)}({k})(-A^2-A^{-2})^{-1}+5f_{i,\hbox{\rm\scriptsize III}}^{(1)}({k})+4A^2f_{i-1,\hbox{\rm\scriptsize III}}^{(1)}({k})(-A^2-A^{-2})\\
      \hskip 19mm +A^4f_{i-2,\hbox{\rm\scriptsize III}}^{(1)}({k})(-A^2-A^{-2})^{2}+(4+A^{-4})f_{i,\hbox{\rm\scriptsize V}}^{(1)}({k})\\
      \hskip 19mm +(4A^2+2A^{-2})f_{i-1,\hbox{\rm\scriptsize V}}^{(1)}({k})(-A^2-A^{-2})+(A^4+1)f_{i-2,\hbox{\rm\scriptsize V}}^{(1)}({k})(-A^2-A^{-2})^{2}.\\
\end{array}
\right.
$$
By combining with (14-15), we get for $i\ge 1$
$$\rho_h(f_{i,j}^{(1)}(k+1))\le 8n+14 \mbox{ and }\rho_l(f_{i,j}^{(1)}(k+1))\ge -4n-10 \hskip 10mm (16)$$
where  $j\in \{\hbox{\rm I},\hbox{\rm II},\hbox{\rm III}, \hbox{\rm IV}\}$.
\noindent By applying Lemma $3.2-3$ with a similar way in the argument of the proof of (16), we obtain for $i\ge 1$
$$\rho_h(f_{i,j}^{(r)}(k+1))\le 8n+14 \mbox{ and }\rho_l(f_{i,j}^{(r)}(k+1))\ge -4n-10 \hskip 10mm (17)$$
where  $j\in \{\hbox{\rm I},\hbox{\rm II},\hbox{\rm III}, \hbox{\rm IV}\}$ and $2\le r\le 3$.

Thus, it is obvious by combining (16-17) that
for $i\ge 1$
$$\rho_h(f_{i}(k+1))\le 8(k+1)+14 \mbox{ and }\rho_l(f_{i}(k+1))\ge -4(k+1)-10.$$
Hence, the conclusion is implied. \hskip 15mm $\Box$
\vskip 3mm
In 1987, Kauffman, Thistlethwaite and Murasugi independently  proved the following result.
\vskip 3mm
\noindent{\bf Lemma $3.5$ (%THEOREM 1 of
\cite{Th87,Ka87,Mu87})} {\it If $L$ is connected, irreducible, alternating link, then the breadth of $V_L(t)$ is precisely $m$.}

\vskip 3mm
\noindent{\bf Proof of Theorem $1.2$.} Since $$<KV_n>=\sum\limits_{i\ge 1}f_{i}(n),$$
it is obvious by Lemma $3.4$ that
for $i\ge 1$
$$\rho_h(<KV_n>)\le 8n+14 \mbox{ and }\rho_l(<KV_n>)\ge -4n-10.$$
Thus,
$$\rho_h(V_{KV_n}(t))\le \dfrac{-3\omega(KV_n)+8n+14}{4} \mbox{ and }\rho_l(V_{KV_n}(t))\ge \dfrac{-3\omega(KV_n)-4n-10}{4}.$$
Then $$br(V_{KV_n}(t))=\rho_h(V_{KV_n}(t))-\rho_l(V_{KV_n}(t))\le 3n+6. \hskip 10mm (18)$$
Because $KV_n$ is connected and irreducible with $4n+6$ crossings,
the result is implied by Lemma 3.5 and the inequality (18). \hskip 10mm $\Box$
\vskip 3mm

\vskip 3mm
\noindent{\bf Proof Theorem $1.3$.} Since $\omega(KV_n)=-4n-6$ and $p_{i}(n)=\sum\limits_{k=1}^3p_{i}^{k}(n)$,
the result is implied by applying a similar way in the argument of the proof of Theorem $1.1$. $\Box$
\vskip 5mm
\noindent{\bf $4.$ Kauffman-Jones polynomials for a type of virtual links }

\vskip 5mm

\vskip 12mm
%\vskip 2mm

\setlength{\unitlength}{0.97mm}
\begin{center}
\begin{picture}(100,50)

\qbezier(-3.5,17)(2,10)(12,16)
\qbezier(5.5,18)(11,19)(12,16)

\qbezier(4,22)(4,15)(4,14)
\qbezier(3.5,18)(-15,18)(-4,22)
\put(-3.5,19){\line(0,1){5}}

\qbezier(-3.5,22)(23,23)(4.5,26)
\qbezier(4,23)(4,27)(4,28)
\put(3.5,26){\line(-1,0){1}}

\put(0,28.5){\circle*{0.7}}
\put(0,26.5){\circle*{0.7}}
\put(0,24.5){\circle*{0.7}}

\qbezier(4,32)(4,30)(4,29.5)
\qbezier(-1.5,28)(-15,28)(-4,32)
\put(-3.5,29){\line(0,1){7}}

\put(-3.5,26){\line(0,1){2}}

\qbezier(-3.5,32)(23,33)(4.5,36)

\qbezier(4,40)(4,36)(4,33)
\qbezier(3.5,36)(-15,36)(-4,40)
\put(-3.5,37){\line(0,1){7}}

\qbezier(-3.5,40)(23,41)(4.5,44)
\qbezier(4,41)(4,45)(4,46)

\qbezier(3.5,44)(-22,48)(-13,20)
\qbezier(-13,20)(-11,10)(3,12)

\qbezier(-3.5,45)(-2,48)(4,46)

\qbezier(46.5,17)(52,10)(62,16)
\qbezier(55.5,18)(61,19)(62,16)

\qbezier(54,22)(54,15)(54,14)
\qbezier(53.5,18)(35,18)(46,22)
\put(46.5,19){\line(0,1){5}}

\qbezier(46.5,22)(73,23)(54.5,26)
\qbezier(54,23)(54,27)(54,28)
\put(53.5,26){\line(-1,0){1}}

\put(50,28.5){\circle*{0.7}}
\put(50,26.5){\circle*{0.7}}
\put(50,24.5){\circle*{0.7}}

\qbezier(54,32)(54,30)(54,29.5)
\qbezier(48.5,28)(35,28)(46,32)
\put(46.5,29){\line(0,1){7}}

\put(46.5,26){\line(0,1){2}}

\qbezier(46.5,32)(73,33)(54.5,36)

\qbezier(54,40)(54,36)(54,33)
\qbezier(53.5,36)(35,36)(46,40)
\put(46.5,37){\line(0,1){7}}

\qbezier(46.5,40)(73,41)(54.5,44)
\qbezier(54,41)(54,45)(54,46)

\qbezier(53.5,44)(28,48)(37,20)
\qbezier(37,20)(39,10)(53,12)

\qbezier(46.5,45)(48,48)(54,46)

\put(36.6,40){\line(0,1){2}}
\put(36.6,40){\line(1,0){2}}

\qbezier(96.5,17)(102,10)(112,16)
\qbezier(105.5,18)(111,19)(112,16)

\qbezier(104,22)(104,15)(104,14)
\qbezier(103.5,18)(85,18)(96,22)
\put(96.5,19){\line(0,1){6}}

\qbezier(96.5,22)(123,23)(104.5,26)
\qbezier(104,23)(104,27)(104,32)
%\put(103.5,26){\line(-1,0){1}}

\qbezier(103.5,26)(85,25)(96,32)

\put(96.5,27){\line(0,1){8}}

\qbezier(96.5,32)(123,33)(104.5,36)

\qbezier(104,40)(104,36)(104,33)
\qbezier(103.5,36)(85,36)(96,40)
\put(96.5,37){\line(0,1){7}}

\qbezier(96.5,40)(123,41)(104.5,44)
\qbezier(104,41)(104,45)(104,46)

\qbezier(103.5,44)(78,48)(87,20)
\qbezier(87,20)(89,10)(103,12)

\qbezier(96.5,45)(98,48)(104,46)

\put(4,13){\circle{1.7}}
\put(54,13){\circle{1.7}}
\put(104,13){\circle{1.7}}

\begin{footnotesize}

\put(-7,45){{$x_{1}$}}
\put(-7,41){{$z_{1}$}}
\put(-7,37){{$z_{2}$}}
\put(-7,32.5){{$z_{3}$}}
\put(-7,26.5){{$z_{4}$}}
\put(4.5,45){{$y_{1}$}}
\put(4.5,41.5){{$y_{2}$}}
\put(4.5,37){{$y_{3}$}}
\put(4.5,30.5){{$y_{4}$}}
\put(4.5,27){{$y_{2n-1}$}}
\put(4.5,23.5){{$y_{2n}$}}
\put(-12,23){{$z_{2n-1}$}}
\put(-8,17){{$z_{2n}$}}
\put(4,11){{$x_{3}$}}
\put(4.5,19){{$x_{2}$}}

\put(43,45){{$x_{1}$}}
\put(43,41){{$z_{1}$}}
\put(43,37){{$z_{2}$}}
\put(43,32.5){{$z_{3}$}}
\put(43,26.5){{$z_{4}$}}
\put(54.5,45){{$y_{1}$}}
\put(54.5,41.5){{$y_{2}$}}
\put(54.5,37){{$y_{3}$}}
\put(54.5,30.5){{$y_{4}$}}
\put(54.5,27){{$y_{2n-1}$}}
\put(54.5,23.5){{$y_{2n}$}}
\put(38,23){{$z_{2n-1}$}}
\put(42,17){{$z_{2n}$}}
\put(54,11){{$x_{3}$}}
\put(54.5,19){{$x_{2}$}}

\put(93,45){{$x_{1}$}}
\put(93,41){{$z_{1}$}}
\put(93,37){{$z_{2}$}}
\put(93,32.5){{$z_{3}$}}
\put(93,27){{$z_{4}$}}
\put(104.5,45){{$y_{1}$}}
\put(104.5,41.5){{$y_{2}$}}
\put(104.5,37){{$y_{3}$}}
\put(104.5,30.5){{$y_{4}$}}
\put(104.5,27){{$y_{5}$}}
\put(104.5,23.5){{$y_{6}$}}
\put(93,23){{$z_{5}$}}
\put(93,16.5){{$z_{6}$}}
\put(104,11){{$x_{3}$}}
\put(104.5,19){{$x_{2}$}}

\put(-5,4){{(a) $RT'_{n}$ }}
\put(40,4){{(b)  oriented $RT'_{n}$ }}
\put(95,4){{(c) $RT'_{3}$ }}
\put(20,-5){{Fig.4: A type of virtual knots $RT'_{n}$ }}
\end{footnotesize}
\end{picture}
\end{center}
\vskip 3mm

In this section, we calculate Kauffman-Jones polynomials of the infinite family of virtual knots $RT'_n$ for $n\ge 1$.
Set $S(RT'_n)$ to be the set of all of states of $RT'_n$ for $n\ge 0$. Denote the following set.

$$
\begin{array}{ll}
S_{\hbox{\rm\scriptsize I}}(RT'_n)=\{s\in S(RT'_n),s_{x_2}=A|\mbox{ either }\exists 1\le k\le n \mbox{ such that }s_{y_{2k-1}}\cdot s_{y_{2k}}\ne A^2, s_{y_i}=A,\\
 \hskip 18mm s_{z_{2k-1}}=s_{z_{2k}}=s_{z_i}=A^{-}\mbox{ for }2k+1\le i\le 2n\mbox{ or }s_{x_1}=s_{z_i}=A^-, s_{y_i}=A \mbox{ for }1\le i\le 2n \},
\end{array}
$$
$
\begin{array}{ll}
S_{\hbox{\rm\scriptsize II}}(RT'_n)=\{s\in S(RT'_n),s_{x_2}=A|\mbox{ either }\exists 1\le k\le n \mbox{ such that }s_{z_{2k-1}}\cdot s_{z_{2k}}\ne A^{-2},  s_{y_i}=A, \\
\hskip 18mm s_{z_i}=A^{-}
\mbox{ for }2k+1\le i\le 2n\mbox{ or }s_{z_i}=A^-, s_{x_1}=s_{y_i}=A \mbox{ for }1\le i\le 2n \},
\end{array}
$
$$
\begin{array}{ll}
S_{\hbox{\rm\scriptsize III}}(RT'_n)=\{s\in S(RT'_n),s_{x_2}=A^-|\mbox{ either }\exists 1\le k\le n \mbox{ such that }s_{y_{2k-1}}\cdot s_{y_{2k}}\ne A^2, s_{y_i}=A,\\
 \hskip 18mm s_{z_{2k-1}}=s_{z_{2k}}=s_{z_i}=A^{-}\mbox{ for }2k+1\le i\le 2n\mbox{ or }s_{x_1}=s_{z_i}=A^-, s_{y_i}=A \mbox{ for }1\le i\le 2n\}
\end{array}
$$
$
\begin{array}{ll}
S_{\hbox{\rm\scriptsize IV}}(RT'_n)=\{s\in S(RT'_n),s_{x_2}=A^-|\mbox{ either }\exists 1\le k\le n \mbox{ such that }s_{z_{2k-1}}\cdot s_{z_{2k}}\ne A^{-2},  s_{y_i}=A, \\
\hskip 18mm s_{z_i}=A^{-}
\mbox{ for }2k+1\le i\le 2n\mbox{ or }s_{x_1}=s_{y_i}=A, s_{z_i}=A^-  \mbox{ for }1\le i\le 2n \}.
\end{array}
$

\vskip 3mm

\noindent{\bf Lemma $4.1$} Let $p_{1,\hbox{\rm\scriptsize I}}(RT'_0)=1$, $p_{1,\hbox{\rm\scriptsize II}}(RT'_0)=A^2$, $p_{2,\hbox{\rm\scriptsize III}}(RT'_0)=A^{-2}$, $p_{1,\hbox{\rm\scriptsize IV}}(RT'_0)=1$. Set $p_{i,\hbox{\rm\scriptsize I}}(RT'_n)=\sum\limits_{s\in S_{\hbox{\rm\scriptsize I}}(RT'_n),l(s)=i}A^{a(s)-b(s)}$,
$p_{i,\hbox{\rm\scriptsize II}}(RT'_n)=\sum\limits_{s\in S_{\hbox{\rm\scriptsize II}}(RT'_n),l(s)=i}A^{a(s)-b(s)}$,
$p_{i,\hbox{\rm\scriptsize III}}(RT'_n)=\sum\limits_{s\in S_{\hbox{\rm\scriptsize III}}(RT'_n),l(s)=i}A^{a(s)-b(s)}$,
$p_{i,\hbox{\rm\scriptsize IV}}(RT'_n)=\sum\limits_{s\in S_{\hbox{\rm\scriptsize IV}}(RT'_n),l(s)=i}A^{a(s)-b(s)}$.
Then for $n\ge 1$
$$
\left\{
\begin{array}{lll}
p_{i,\hbox{\rm\scriptsize I}}(RT'_n)=p_{i,\hbox{\rm\scriptsize
      I}}(RT'_{n-1})+2A^{-2}p_{i-1,\hbox{\rm\scriptsize I}}(RT'_{n-1})+A^{-4}p_{i-2,\hbox{\rm\scriptsize I}}(RT'_{n-1})\\
      \hskip 19mm +2A^{-2}p_{i,\hbox{\rm\scriptsize II}}(RT'_{n-1})+A^{-4}p_{i-1,\hbox{\rm\scriptsize II}}(RT'_{n-1}); \\
p_{i,\hbox{\rm\scriptsize II}}(RT_n)=2A^2p_{i,\hbox{\rm\scriptsize
      I}}(RT'_{n-1})+(A^4+4)p_{i-1,\hbox{\rm\scriptsize I}}(RT'_{n-1})\\
      \hskip 19mm +(2A^2+2A^{-2})p_{i-2,\hbox{\rm\scriptsize I}}(RT'_{n-1})+p_{i-3,\hbox{\rm\scriptsize I}}(RT'_{n-1})+5p_{i,\hbox{\rm\scriptsize II}}(RT'_{n-1})\\
      \hskip 19mm +(4A^2+2A^{-2})p_{i-1,\hbox{\rm\scriptsize II}}(RT'_{n-1})+(A^4+1)p_{i-2,\hbox{\rm\scriptsize II}}(RT'_{n-1}); \\
p_{i,\hbox{\rm\scriptsize III}}(RT_n)=p_{i,\hbox{\rm\scriptsize
      III}}(RT'_{n-1})+2A^{-2}p_{i-1,\hbox{\rm\scriptsize III}}(RT'_{n-1})+A^{-4}p_{i-2,\hbox{\rm\scriptsize III}}(RT'_{n-1})\\
      \hskip 19mm +2A^{-2}p_{i-1,\hbox{\rm\scriptsize IV}}(RT'_{n-1})+A^{-4}p_{i-2,\hbox{\rm\scriptsize IV}}(RT'_{n-1}); \\
p_{i,\hbox{\rm\scriptsize IV}}(RT'_n)=2A^2p_{i+1,\hbox{\rm\scriptsize
      III}}(RT'_{n-1})+(A^4+4)p_{i,\hbox{\rm\scriptsize III}}(RT'_{n-1})\\
      \hskip 19mm +(2A^2+2A^{-2})p_{i-1,\hbox{\rm\scriptsize III}}(RT'_{n-1})+p_{i-2,\hbox{\rm\scriptsize III}}(RT'_{n-1})+5p_{i,\hbox{\rm\scriptsize IV}}(RT'_{n-1})\\
      \hskip 19mm +(4A^2+2A^{-2})p_{i-1,\hbox{\rm\scriptsize IV}}(RT'_{n-1})+(A^4+1)p_{i-2,\hbox{\rm\scriptsize IV}}(RT'_{n-1}).
\end{array}
\right.
$$
%\vskip 3mm

\vskip 3mm

\noindent{\bf Proof of Theorem $1.4$} By a similar way in the argument of the proof of Theorem $1.1$, the result holds.

\hskip 150mm $\Box$

\vskip 5mm
 \hskip 55mm{\bf $5.$ Further study }
\vskip 5mm
%\vskip 2mm
In this section, we introduce general $m$-string alternating links(or virtual links) and $m$-string tangle links(or virtual links) for $m\ge 2$. Several problems are proposed.

   Set $n$ to be a positive integer in this section. Generally, given a link(or virtual link) $L_0$ with $m+1$ parallel edges, denote one of them by $e_0$ and denote others by $e_i=(u_i^{r_i},v_i^{\varepsilon_i})$ in sequence for $1\le i\le m$. Here $r_i\in\{+,-\}$, $r_0=r$. Assume that $e_0$ is the leftmost edge and leave other cases to readers to get in a similar way. Add $2n$ crossings $x_{i,l}$ on $e_i$ in sequence for $1\le i\le m$ and $1\le l\le 2n$ respectively. Let $(u_{i}^{r_i},x_{i,1}^{r})$, $(x_{i,1}^{r},x_{i,2}^{-r})$, $\cdots$, $(x_{i,2n-1}^{r},x_{i,2n}^{-r})$, $(x_{i,2n}^{-r},v_{i}^{\epsilon_{i}})$ be a subdivision of $e_i$ for odd $1\le i\le m$ and let $(u_{i}^{r_i},x_{i,1}^{-r})$, $(x_{i,1}^{-r},x_{i,2}^{r})$, $\cdots$, $(x_{i,2n-1}^{-r},x_{i,2n}^{r})$, $(x_{i,2n}^{r},v_{i}^{\epsilon_{i}})$ be a subdivision of $e_i$ for even $1\le i\le m$. Add edges $(u_{i}^{r_{i}},x_{i,1}^{-r})$ and $(v_{i}^{\varepsilon_{i}},x_{i,2n}^{r})$ for odd $1\le i\le m$, add edges $(u_{i}^{r_{i}},x_{i,1}^{r})$ and $(v_{i}^{\varepsilon_{i}},x_{i,2n}^{-r})$ for even $1\le i\le m$, add edges $(x_{i,l}^{-r},x_{i+1,l}^{r})$, $(x_{i,l}^{-r},x_{i-1,l}^{r})$ for odd $1\le i\le m$ and odd $1\le l\le 2n$,
$(x_{i,l}^{r},x_{i+1,l}^{-r})$, $(x_{i,l}^{r},x_{i-1,l}^{-r})$ for odd $1\le i\le m$ and even $1\le l\le 2n$, and then add edges $(x_{m,l}^r,x_{m,l+1}^{-r})$ for even $m$ and odd $1\le l\le 2n$.
A link (or a virtual link) $AL_n$ is constructed which is  called an {\it $m$-string alternating} link(or virtual link). Here,
$x_{m+1,l}^{r}=x_{m,l+1}^{r}$ and $x_{0,l}^{r}=x_{1,l-1}^{r}$ for odd $1\le l\le 2n$, $x_{m+1,l}^{-r}=x_{m,l-1}^{-r}$ and $x_{0,l}^{-r}=x_{1,l+1}^{-r}$ for even $1\le l\le 2n$, $x_{0,1}^r=x_{1,0}^r=u_0^r$, $x_{0,2n+1}^{-r}=v_0^{\varepsilon_0}$. An example is shown in Fig.5 (b)  for $m=3$ and $n=2$.
%%%%%%%%%%%%%%%%%%%%%%%%%%%%%%%%%%%%%%%%%%%%%%%%%%%%%%%%%%%%%%%%%%%%%%%%%%%%%%%%%%%%%%%
\vskip 12mm
%\vskip 2mm

\setlength{\unitlength}{0.97mm}
\begin{center}
\begin{picture}(100,50)
\put(9,14){\line(0,1){33}}
\put(-1,14){\line(0,1){33}}
\put(-11,14){\line(0,1){33}}
\put(-21,14){\line(0,1){33}}
\put(59,14){\line(0,1){6}}
\put(59,21){\line(0,1){15}}
\put(59,37){\line(0,1){11}}

\put(49,14){\line(0,1){12}}
\put(49,28){\line(0,1){14}}
\put(49,44){\line(0,1){4}}

\put(39,14){\line(0,1){4}}
\put(39,20){\line(0,1){14}}
\put(39,35){\line(0,1){13}}

%\put(29,17){\line(0,1){30}}
\put(29,43){\line(1,0){9}}
\put(40,43){\line(1,0){18}}
%\put(50,43){\line(1,0){8}}

\qbezier(60,43)(70,38)(50,35)
\put(48,35){\line(-1,0){8}}
\qbezier(40,35)(33,30)(38,27)
\put(40,27){\line(1,0){10}}

\put(50,27){\line(1,0){8}}
\qbezier(60,27)(70,22)(50,19)
\put(48,19){\line(-1,0){18}}
%\put(38,19){\line(-1,0){8}}

\put(109,14){\line(0,1){6}}
\put(109,21){\line(0,1){15}}
\put(109,37){\line(0,1){11}}

%\put(99,14){\line(0,1){34}}
\put(99,14){\line(0,1){4}}
\put(99,19){\line(0,1){14}}
\put(99,35){\line(0,1){13}}

%\put(89,17){\line(0,1){30}}
\put(89,14){\line(0,1){4}}
\put(89,19){\line(0,1){14}}
\put(89,35){\line(0,1){13}}

%\put(29,17){\line(0,1){30}}
\put(79,43){\line(1,0){9}}
\put(90,43){\line(1,0){8}}
\put(100,43){\line(1,0){8}}

\qbezier(110,43)(120,38)(100,35)
\put(100,35){\line(-1,0){10}}
%\qbezier(89,35)(78,30)(98,27)
\qbezier(90,35)(83,30)(88,27)
\put(90,27){\line(1,0){8}}
\put(100,27){\line(1,0){8}}
\qbezier(110,27)(120,22)(100,19)
%\put(98,19){\line(-1,0){8}}
\put(100,19){\line(-1,0){20}}

\put(-23,10){$v_0^{\varepsilon_0}$}
\put(-13,10){$v_1^{\varepsilon_1}$}
\put(-3,10){$v_2^{\varepsilon_2}$}
\put(7,10){$v_3^{\varepsilon_3}$}
\put(-23,50){$u_0^{r}$}
\put(-13,50){$u_1^{r_1}$}
\put(-3,50){$u_2^{r_2}$}
\put(7,50){$u_3^{r_3}$}

\put(-11,4){(a) $L_0$}

\put(27,15){$v_0^{\varepsilon_0}$}
\put(37,10){$v_1^{\varepsilon_1}$}
\put(47,10){$v_2^{\varepsilon_2}$}
\put(57,10){$v_3^{\varepsilon_3}$}
\put(27,45){$u_0^{r}$}
\put(37,50){$u_1^{r_1}$}
\put(47,50){$u_2^{r_2}$}
\put(57,50){$u_3^{r_3}$}

\put(32,40){$x_{1,1}$}
\put(42,40){$x_{1,2}$}
\put(52,40){$x_{1,3}$}

\put(31,33){$x_{2,1}$}
\put(49,32){$x_{2,2}$}
\put(59,33){$x_{2,3}$}

\put(31,26){$x_{3,1}$}
\put(43,28.5){$x_{3,2}$}
\put(59,28){$x_{3,3}$}

\put(31,21){$x_{4,1}$}
\put(43,21){$x_{4,2}$}
\put(60,19){$x_{4,3}$}

\put(39,4){(b) $AL_4$}

\put(77,15){$v_0^{\varepsilon_0}$}
\put(87,10){$v_1^{\varepsilon_1}$}
\put(97,10){$v_2^{\varepsilon_2}$}
\put(107,10){$v_3^{\varepsilon_3}$}
\put(77,45){$u_0^{r}$}
\put(87,50){$u_1^{r_1}$}
\put(97,50){$u_2^{r_2}$}
\put(107,50){$u_3^{r_3}$}

\put(82,40){$x_{1,1}$}
\put(92,40){$x_{1,2}$}
\put(102,40){$x_{1,3}$}

\put(81,33){$x_{2,1}$}
\put(99,32){$x_{2,2}$}
\put(109,33){$x_{2,3}$}

\put(81,26){$x_{3,1}$}
\put(93,28.5){$x_{3,2}$}
\put(109,28){$x_{3,3}$}

\put(81,21){$x_{4,1}$}
\put(93,21){$x_{4,2}$}
\put(110,19){$x_{4,3}$}

\put(89,4){(c) $TL_4$}
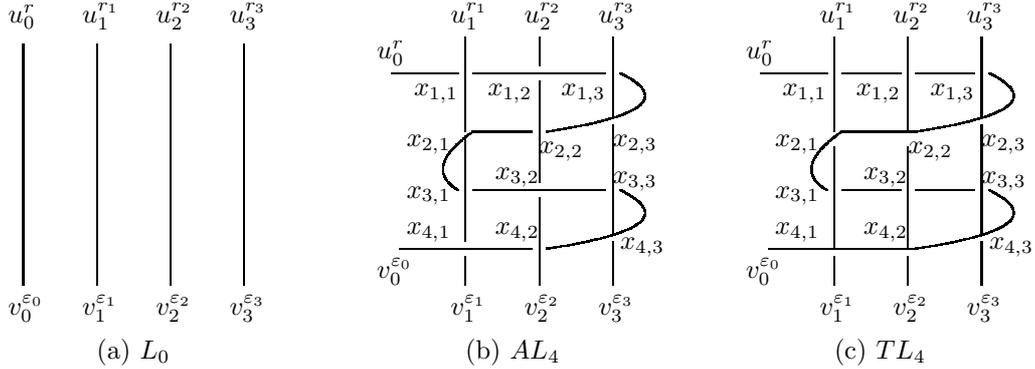
\put(29,-2){Fig.$5$:$L_0$, $AL_4$ and $TL_4$}
\end{picture}
\end{center}

%%%%%%%%%%%%%%%%%%%%%%%%%%%%%%%%%%%%%%%%%%%%%%%%%%%%%%%%%%%%%%%%%%%%%%%%%%%%%%%%%%%%%%%

\vskip 3mm
Similarly,
let $L_0$ be a link (or virtual link) above. Add $2n$ crossings $x_{i,l}$ on $e_i$ in sequence for $1\le i\le m$ and $1\le l\le 2n$ respectively. Let $(u_{i}^{r_i},x_{i,1}^{r})$, $(x_{i,1}^{r},x_{i,2}^{-r})$, $\cdots$, $(x_{i,2n-1}^{r},x_{i,2n}^{-r})$, $(x_{i,2n}^{-r},v_{i}^{\epsilon_{i}})$ be a subdivision of $e_i$.  Add edges $(u_{i}^{r_{i}},x_{i,1}^{r})$ and $(v_{i}^{\varepsilon_{i}},x_{i,2n}^{-r})$ for $1\le i\le l$,
add edges $(x_{i,l}^{-r},x_{i+1,l}^{-r})$, $(x_{i,l}^{-r},x_{i-1,l}^{-r})$ for odd $1\le i\le m$ and odd $1\le l\le 2n$,
$(x_{i,l}^{r},x_{i+1,l}^{r})$, $(x_{i,l}^{r},x_{i-1,l}^{r})$ for odd $1\le i\le m$ and even $1\le l\le 2n$, and then add edges $(x_{m,l}^{-r},x_{m,l+1}^{r})$ for even $m$ and odd $1\le l\le 2n$.
A link (or virtual link) $TL_n$ is constructed which is  called an {\it $m$-string tangle} link (or virtual link). Here,
$x_{m+1,l}^{-r}=x_{m,l+1}^{r}$ and $x_{0,l}^{-r}=x_{1,l-1}^{r}$ for odd $1\le l\le 2n$, $x_{m+1,l}^{r}=x_{m,l-1}^{-r}$ and $x_{0,l}^{r}=x_{1,l+1}^{-r}$ for even $1\le l\le 2n$, $x_{0,1}^{-r}=x_{1,0}^{-r}=u_0^r$, $x_{0,2n+1}^{r}=v_0^{\varepsilon_0}$. An example is also shown in Fig.5 (c) for $m=3$ and $n=2$.

\vskip 3mm
\noindent{\bf Problem $5.1.$ }{\it Given a link $L_0$,  let $AL_n$ is an $m$-string alternating link constructed from $L_0$ for $m\ge 3$. Determine $V_{AL_n}(t)$.}

\vskip 3mm
\noindent{\bf Problem $5.2.$ }{\it Given a link $L_0$,  let $TL_n$ is an $m$-string tangle link constructed from $L_0$ for $m\ge 3$. Determine $V_{TL_n}(t)$.}

\vskip 3mm
\noindent{\bf Conjecture $5.3.$ }{\it Suppose that $L_0$ is connected and  irreducible and that $TL_n$ is an $m$-string tangle link constructed from $L_0$ for $m\ge 2$. If $L_0$ is non-alternating, then $TL_n$ is also non-alternating. }

\vskip 3mm
\noindent{\bf Conjecture $5.4.$ }{\it Suppose that a link $L_0$ is prime  and that $L_n$ is an $m$-string alternating (or tangle) link constructed from $L_0$ for $m\ge 2$. Then $L_n$ is prime. }

\vskip 5mm

%\vskip 5mm

\end{document}